\documentclass[11pt]{article}
\usepackage{amsfonts}
\usepackage[mathscr]{eucal}
\usepackage{mathrsfs}
\usepackage{amsmath}
\usepackage{amssymb}
\usepackage{amsmath}
\usepackage{amsthm}
\usepackage{mathrsfs}
\usepackage{graphicx}
\usepackage{comment}
\usepackage{url,hyperref,lineno,microtype,subcaption}
\usepackage[onehalfspacing]{setspace}

\newtheorem{theorem}{Theorem}[section]
\newtheorem{corollary}{Corollary}

\newtheorem{definition}[theorem]{Definition}
\newtheorem{remark}{Remark}

\newcommand{\Ga}{{\Gamma}}

\newcommand{\Om}{{\Omega}}

\newcommand{\be}{{\beta}}
\newcommand{\ga}{{\gamma}}

\newcommand{\ep}{\varepsilon}

\newcommand{\la}{\lambda}
\newcommand{\si}{{\sigma}}

\newcommand{\vph}{{\varphi}}
\newcommand{\om}{{\omega}}

\newcommand{\R}{{\mathbb R}}

\newcommand{\cE}{\mathcal{E}}

\newcommand{\cK}{\mathcal{K}}

\newcommand{\cW}{\mathcal{W}}

\newcommand{\pd}{\partial}
\newcommand{\intl}{\int\limits}
\newcommand{\arr}{\rightarrow}

\newcommand{\oPhi}{\overline{\Phi}}

\newcommand{\oom}{\overline{\om}}
\newcommand{\ovph}{\overline{\varphi}}
\newcommand{\ou}{\overline{u}}
\newcommand{\ov}{\overline{v}}
\newcommand{\ow}{\overline{w}}

\title {On the asymptotic properties of solutions to a nonlinear transmission problem for a Bresse beam with thermal damping.}

\author{Tamara Fastovska\,$^{1,3,*}$ and Dirk Langemann\,$^{2}$ } 
\date{} 
\begin{document}
	
	\onecolumn

	

	\maketitle
	
	{\bf Abstract.}
	A  nonlinear transmission problem for an arch beam, which consists of two parts with different material properties is considered. One of the parts is subjected to thermal damping while another one is undamped. The thermal damping affects only the shear angle and the longitudinal equation within the damped part, which can have any positive length. We prove the well-posedness of the system in the energy space and establish the existence of a compact global attractor under certain conditions imposed solely on the coefficients of the damped part. Additionally, we study singular limits of the system as the beam’s curvature tends to zero or in case both the curvature vanishes and the shear moduli tend to infinity. Numerical simulations are also carried out to model these limiting behaviors.

{\bf 2020 Mathematics Subject Classification:} Primary: 35B41, 35B25; Secondary:74K10. 

{ {\bf Keywords:} thermoelastic  Bresse beam,  transmission problem, global attractor, singular limit.  }

\section{Introduction}

\begin{figure}[h]
	\centering
	\includegraphics[width=0.5\textwidth]{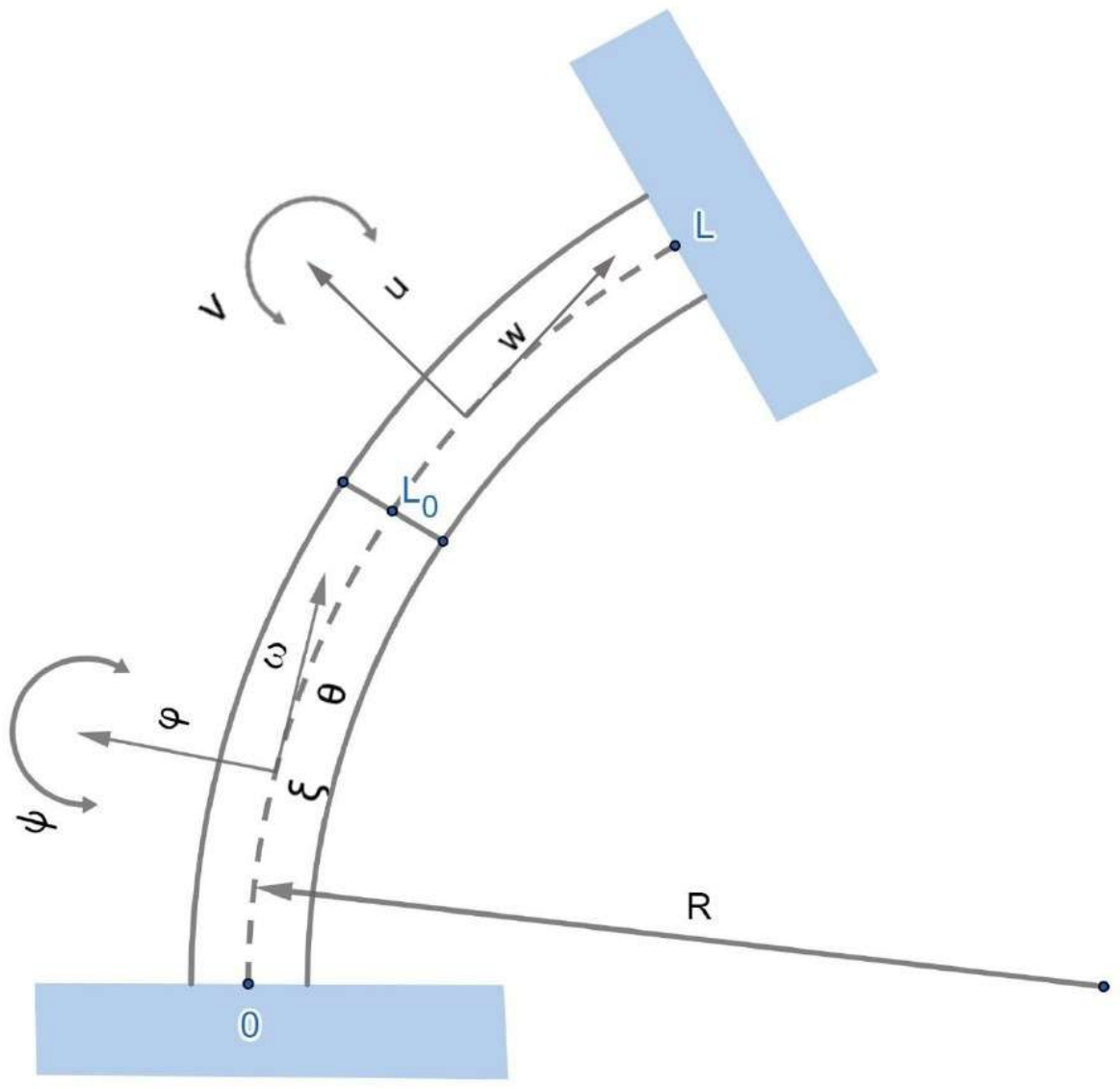}
	\caption{Composite Bresse beam with temperature dissipation.} \label{FigBeam}
\end{figure}

The Bresse model describes the dynamics of an initially curved beam while accounting for shear deformation. The original formulation of the model for a homogeneous Bresse beam was introduced in \cite{Bre1859}. In this work, we study a variant of the model presented in \cite{You2022}, incorporating nonlinear effects as outlined in \cite[Ch. 3]{LagLeu1994}.  Let the whole beam occupy a part of a circle of length $L$ and have the curvature $l=R^{-1}$. Treating the beam as a one-dimensional object, we define the spatial variable $x$ along its length, ranging over the interval $(0,L)$.  The parts of the beam occupying the intervals $(0,L_0)$ and  $(L_0,L)$ consist of  materials with different material properties. The part lying in the interval $(0,L_0)$ is partially subjected to a thermal damping  (Figure \ref{FigBeam}).

We denote by $\vph$, $\psi$, and $\om$  the transversal displacement, the shear angle variation, and  the longitudinal displacement of the left part of the beam lying in $(0,L_0)$. Analogously, we denote by  $u$, $v$, and $w$  the transversal displacement,  the shear angle variation,  and the longitudinal displacement of the right part of the beam occupying the interval $(L_0,L)$. The unknown functions $\xi$ and $\theta$ represent  the temperature deviations
from a reference temperature along the longitudinal and vertical directions \cite{You2022}.    We assume the beam is fixed at both ends. For the temperature variables, Dirichlet boundary conditions are imposed at both endpoints of the interval  $(0,L_0)$.
Nonlinear oscillations of the composite beam can be described by the following  system of equations
\begin{align}
	& \rho_1\vph_{tt}-k_1(\vph_x+\psi+l\om)_x - l\si(\om_x-l\vph)-\frac{\sigma l}{2}\psi^2+l\alpha_1\xi=p_1(x,t),\;x\in (0,L_0), t>0,\label{Eq1}\\
	& \be_1\psi_{tt} -\nu_1 \psi_{xx} +k_1(\vph_x+\psi+l\om) +\alpha_2\theta_x+\sigma\psi(\omega_x-l\vph)-\alpha_1\psi\xi+\frac{\sigma}{2}\psi^3=r_1(x,t),\label{Eq2}\\
	& \rho_1\om_{tt}- \si(\om_x-l\vph) _x+lk_1(\vph_x+\psi+l\om)-\sigma\psi\psi_x+\alpha_1\xi_x=q_1(x,t), \label{Eq3}\\
	&\gamma\xi_t-\mu\xi_{xx} +\alpha_1(\omega_x-l\vph)_t+\alpha_1\psi_t\psi=h_1(x,t),\label{Eq31}\\
	&\delta\theta_t-\lambda\theta_{xx}+\alpha_2\psi_{tx}=h_2(x,t)\label{Eq32}
\end{align}
and
\begin{align}
	& \rho_2u_{tt}-k_2(u_x+v+lw)_x - l\si(w_x-lu) -\frac{\sigma    l}{2}v^2=p_2(x,t), \label{Eq4}\\
	& \be_2v_{tt}-\nu_2 v_{xx}+k_2(u_x+v+lw)+\sigma v(w_x-lu)+\frac{\sigma}{2}v^3 =r_2(x,t),\qquad x\in (L_0,L), t>0,\label{Eq5}\\
	& \rho_2w_{tt}- \si(w_x-lu)_x+lk_2(u_x+v+lw)-\sigma vv_x=q_2(x,t), \label{Eq6}
\end{align}
where $\rho_j,\;\be_j, \; k_j, \; \si,\; \nu_j,\; \alpha_j,\; \delta,\;\gamma,\;\lambda,\;\mu$ are positive physical parameters,  $p_j, \; q_j, \; r_j$ are known external loads, $h, \; g,$ are known external heat sources.
The system is subjected to the Dirichlet boundary conditions at the ends of the beam
\begin{align} 
	\vph(0,t)=u(L,t)=0,  \quad &\psi(0,t)=v(L,t)=0, \quad \om(0,t)=w(L,t)=0,\label{BC}\\
	&\theta(0,t)=\xi(0,t)=0	.\label{BC1}
\end{align}
and transmission conditions at point $L_0$
\begin{align}
	& \vph(L_0,t)=u(L_0,t), \quad  \psi(L_0,t)=v(L_0,t),  \quad  \om(L_0,t)=w(L_0,t), \label{TC1} \\
	& k_1(\vph_x+\psi+l\om)(L_0,t)=k_2(u_x+v+lw)(L_0,t), \\
	& \nu_1 \psi_{x}(L_0,t)= \nu_2 v_{x}(L_0,t),\\
	& \om_x(L_0,t)=w_x(L_0,t),\label{TC44}\\
	&\xi(L_0,t)=\theta(L_0,t)=0.\label{TC4}
\end{align}
The system is also supplemented with the initial conditions
\begin{align}
	&\vph(x,0)=\vph_0(x),\quad \psi(x,0)=\psi_0(x),\quad \om(x,0)=\om_0(x),\quad \xi(x,0)=\xi_0(x),\\
	&\vph_t(x,0)=\vph_1(x),\quad \psi_t(x,0)=\psi_1(x),\quad \om_t(x,0)=\om_1(x),\quad \theta(x,0)=\theta_0(x),
\end{align}
\begin{align}
	&u(x,0)=u_0(x),\quad v(x,0)=v_0(x),\quad w(x,0)=w_0(x),\\
	&u_t(x,0)=u_1(x),\quad v_t(x,0)=v_1(x),\quad w_t(x,0)=w_1(x).\label{IC}
\end{align}
In what follows we use the notations for some physical expressions:
\begin{align*}
	& Q_i(\nu, \zeta, \eta)=k_i(\xi_x+\zeta +l\eta) \mbox{  are  shear forces},\\
	& N(\nu, \zeta, \eta)=\si(\eta_x-l\xi) \mbox{ are axial forces},\\
	& M_i(\nu, \zeta, \eta)=\nu_i\zeta_x \mbox{ are bending moments}
\end{align*}
for damped ($i=1)$ and undamped ($i=2$) parts respectively.
Later we will use them to rewrite the problem in an abstract form.

In this paper, we study the well-posedness and long-time behavior of the system \eqref{Eq1}–\eqref{IC}. Our primary objective is to identify conditions on the system's coefficients under which thermal dissipation, applied only to one part of the beam, is sufficient to ensure the existence of a global attractor, regardless of the length of the damped part. Additionally, we analyze singular limits of the Bresse system as the curvature parameter 
$l$ tends to zero. In this case, solutions converge to those of a nonlinear thermoelastic Timoshenko system, which models the dynamics of a straight, shearable beam. Furthermore, in the double limit where the curvature 
$l$ tends to zero and the shear moduli $k_1, k_2$ tend to infinity, the system converges to a transmission full von Kármán model with rotational inertia, coupled with heat conduction equations. These asymptotic behavior is illustrated through numerical simulations.

The long-time behaviour of linear homogeneous Bresse beams with various boundary conditions and dissipation types was profoundly investigated. 
If damping influences  all three equations for unknown displacements and rotation, it is sufficient  for the  exponential stability of the system without additional assumptions on the parameters of the problems (see, e.g., \cite{AlmSan2010, Mukiawa1, Mukiawa2}).

The situation may change if we have a dissipation of any kind  in one or two equations only.  There are examples  for the Timoshenko beams \cite{MuRa2002} and for the Bresse beams  \cite{Oro2015, LiuRao2009} which show that damping in only one of the equations  does not guarantee the exponential stability of the  system.  The exponential stability is guaranteed only under  additional assumptions on the coefficients of the problem, usually, the equality of the propagation speeds 
\begin{equation*}
	k_1=\sigma, \quad 	\frac{\rho_1}{k_1}=\frac{\beta_1}{\lambda_1}.
\end{equation*}
If these relations do not hold true, only polynomial stability can be shown (\cite{AlaMuAl2011} in case of mechanical dissipation and \cite{Oro2015} in case of thermal dissipation).

In the present paper we study a transmission problem for the Bresse system. Transmission problems for various equation types have already received an ample attention. One can find a number of papers concerning their well-posedness and long-time behaviour of linear transmission problems for elastic bodies with various types of damping. For instance, exponential stability  was proved under certain conditions on parameters for a Kirchhoff-Timoshenko thermoelastic system \cite{Fast2013}, a wave equation with thermal damping \cite{FaLuMu2003}, Kirchhoff type plates and beams \cite{MuPo2004,MuPo2011}. Polynomial stability was shown for a thermoelasic Bresse systems with dynamic transmission conditions in \cite{You2022}. 

For nonlinear transmission problems there are much less results on the qualitative behaviour of solutions, in particular, on the existence of attractors (see, e.g., \cite{Pot2012} for a  thermoelastic plate with  Berger's nonlinearity,  \cite{Fast2023} for the full von Kármán beam with structural damping and \cite{Fast2024} for a semilinear Bresse model with structural damping). The combination of transmission boundary conditions and a nonlinearity  creates substantial technical difficulties.

To the best of our knowledge, there are only a few  papers  on the long-time behaviour of transmission problems for the Bresse system.  The beam in \cite{You2022} consists of a thermoelastic (damped) and elastic (undamped) parts, both purely linear with dynamic transmission conditions.  It was shown that the corresponding semigroup is not exponentially stable for any set of parameters, but only polynomially  stable. In \cite{Fast2024} a semilinear transmission Bresse problem with mechanical damping only in the equation for the shear angle for the damped part is investigated. The existence of a compact global attractor  under  restrictions on the coefficients in the damped part is established. 

In the present work we consider a transmission problem for the nonlinear Bresse system with thermal dissipation acting on one of the parts of the beam. Thermal damping is weaker than mechanical in principal, moreover, and in the specific case being analyzed, its effect on transverse displacement is even smaller. Therefore, the question to be answered is whether the damping acting on the longitudinal displacement and the shear angle is sufficient for the dissipation of the transversal displacement and the whole undamped part of the beam. In this case the Dirichlet boundary conditions for the temperature variations create significant technical obstacles in the proof of the asymptotic smoothness property. We may observe an interesting effect of boundary conditions for the heat equation in the transmission thermoelastic problems. It turns out that in case of the Neumann or Robin boundary conditions at least at one of the boundary points (see e.g. \cite{FaLuMu2003,MuPo2004}) the proof of asymptotic smoothness in the nonlinear case (or exponential stability in linear case) is easier than in case  of the Dirichlet conditions at both ends. In the latter case (see, \cite{MuPo2011}) one has to apply more sophisticated multipliers. However, for the system under our consideration even this approach does not work, since we have an undamped transversal component in the thermoelastic part of the system.

The paper is organized as follows. In Section 2 we represent functional spaces and pose the problem in an abstract form. In Section 3 we prove that the problem is well-posed and possesses strong solutions provided nonlinearities and initial data are smooth enough.  Section 4 is devoted  to the main result on the existence of a compact attractor. The  combination of weak and localized dissipation with the presence of nonlinearity  prevents us from proving the existence of an absorbing ball explicitly, thus we show that the dynamical system generated is  gradient  and asymptotically smooth under certain assumptions on coefficients of the problem.  In Section 5 we show that solutions to \eqref{Eq1}-\eqref{IC} tend to solutions to a thermoelastic transmission problem for a nonlinear Timoshenko beam when $l\arr 0$ and to solutions to  transmission problem for a thermoelastic full von Karman beam with rotational inertia when $l\arr 0$ and $k_i \arr\infty$ as well as perform numerical modelling of these singular limits.

\section{Preliminaries and abstract formulation}
Let us denote
\begin{equation*}
	\Phi^1=(\vph, \psi, \om), \quad\Phi^2=(u,v,w),\quad \Theta=(\xi, \theta), \quad \Phi=(\Phi^1,\Phi^2).
\end{equation*}
Analogously,
\begin{align*}
	& P_j=(p_j,q_j, r_j): [(0,L)\times \R_+]\arr\R^3,\quad \EuScript P=(P_1,P_2),\quad G=(h_1, h_2): [(0,L_0)\times \R_+]\arr\R^2  , \\
	& S=diag\{\gamma,\delta\}, \quad R_j=diag\{\rho_j,\be_j,\rho_j\},\quad R=diag\{R_1, R_2\}.
\end{align*}
where $j=1,2$. \\
Throughout  the paper we use the notation $||\cdot||$ for the $L^2$-norm of a function and $(\cdot,\cdot)$ for the $L^2$-inner product. In these notations we skip  the domain, on which functions are defined. We adopt the notation $||\cdot||_{L^2(\Omega)}$ only when the domain is not evident.  We also use the same notations  $||\cdot||$ and  $(\cdot,\cdot)$ for $[L^2(\Om)]^3$. We also write $||\cdot||_\Omega$ to highlight the domain of integration.\\
To write our problem in an abstract form we introduce the following spaces. For the velocities of the displacements we use the space
\begin{equation*}
	H_v=\{\Phi=(\Phi^1,\Phi^2):\; \Phi^1\in [L^2(0,L_0)]^3, \;\Phi^2\in [L^2(L_0,L)]^3 \}
\end{equation*}
with the norm
\begin{equation*}
	||\Phi||^2_{H_v}=||\Phi||^2_v=\sum_{j=1}^{2}||\sqrt{R_j}\Phi^j||^2
\end{equation*}
which is equivalent to the standard $L^2$-norm.\\
For the beam displacements we use the space
\begin{align*}
	H_d=\left\{\Phi\in H_v:\; \right.&\Phi^1\in [H^1(0,L_0)]^3, \;\Phi^2\in [H^1(L_0,L)]^3,  \\
	&\left.\Phi^1(0,t)=\Phi^2(L,t)=0,\; \Phi^1(L_0,t)=\Phi^2(L_0,t) \right\}
\end{align*}
with the norm
\begin{equation*}
	||\Phi||^2_{H_d}=||\Phi||^2_d=\sum_{j=1}^{2}\left(||Q_j(\Phi^j)||^2+||N(\Phi^j)||^2+||M_j(\Phi^j)||^2 \right).
\end{equation*}
This norm is equivalent to the standard $H^1$-norm. Moreover, the equivalence constants can be chosen independent of $l$ for $l$ small enough (see \cite{MaMo2017}, Remark 2.1). If we  set
\begin{equation*}
	\Psi(x)=\left\{
	\begin{aligned}
		&\Phi^1(x), \quad &x\in (0,L_0),\\
		&\Phi^2(x), \quad &x\in [L_0,L)
	\end{aligned}
	\right.
\end{equation*}
we see that there is isomorphism between $H_d$ and $[H^1_0(0,L)]^3$.
For the thermal components we utilize
\begin{equation*}
	H_\theta=\{\Theta=(\xi,\theta):\; \xi, \theta\in L^2(0,L_0)\}=[L^2(0,L_0)]^2.
\end{equation*}
with the norm
\begin{equation*}
	||\Theta||^2_{H_\theta}=||\Theta||^2_\theta=\gamma||\xi||^2+\delta ||\theta||^2
\end{equation*}
and 
\begin{equation*}
	\tilde H_\theta=\{\Theta=(\xi,\theta):\; \xi, \theta\in H_0^1(0,L_0)\}=[H_0^1(0,L_0)]^2.
\end{equation*}
with the norm
\begin{equation*}
	||\Theta||^2_{\tilde H_\theta}=\mu||\xi_x||^2+\lambda ||\theta_x||^2.
\end{equation*}
For the initial conditions we use the notations
\begin{align*}
	&\Phi_0(x)=(\vph_0(x),\psi_0(x),\omega_0(x), u_0(x), v_0(x), w_0(x)), \\
	&\Phi_1(x)=(\vph_1(x),\psi_1(x),\omega_1(x), u_1(x), v_1(x), w_1(x)),\\
	&\Theta_0=(\xi_0(x), \theta_0(x)).
\end{align*}
We consider an operator  $B:D(B)\subset H_v\arr H_v$
\begin{equation}\label{oper}
	B\Phi=\left(
	\begin{aligned}
		& -\pd_x Q_1(\Phi^1)-lN(\Phi^1)\\
		& -\pd_x M_1(\Phi^1)+Q_1(\Phi^1)\\
		& -\pd_x N(\Phi^1)+lQ_1(\Phi^1) \\
		& -\pd_x Q_2(\Phi^2)-lN(\Phi^2)\\
		& -\pd_x M_2(\Phi^2)+Q_2(\Phi^2) \\
		& -\pd_x N(\Phi^2)+lQ_2(\Phi^2)
	\end{aligned}
	\right)
\end{equation}
with the domain
\begin{multline*}
	D(B)=\left\{\Phi\in H_d:\right. \Phi^1\in H^2(0,L_0), \; \Phi^2\in H^2(L_0,L), \;
	Q_1(\Phi^1(L_0,t))=Q_2(\Phi^2(L_0,t)), \\
	N(\Phi^1(L_0,t))=N(\Phi^2(L_0,t)),
	\left.M_1(\Phi^1(L_0,t))=M_2(\Phi^2(L_0,t))\right.\}.
\end{multline*}
Arguing analogously to  Lemma 2.1 from \cite{Fast2024} one can prove that the operator $A$ is positive  and self-adjoint. Moreover,
$(B^{1/2}\Phi, B^{1/2}\tilde \Phi)=\EuScript B(\Phi, \tilde \Phi)$
\begin{equation*} \label{AForm}
	\begin{aligned}
		\EuScript B(\Phi, \tilde \Phi)	=
		& \frac{1}{k_1} (Q_1(\Phi^1),Q_1(\tilde\Phi^1))_{[0,L_0]} + \frac{1}{\si} (N(\Phi^1),N(\tilde\Phi^1))_{[0,L_0]} + \frac{1}{\la_1} (M_1(\Phi^1),M_1(\tilde\Phi^1))_{[0,L_0]} + \\
		& \frac{1}{k_2} (Q_2(\Phi^2),Q_2(\tilde\Phi^2))_{[L_0,L]} + \frac{1}{\si} (N(\Phi^2),N(\tilde\Phi^2))_{[L_0,L]} + \frac{1}{\la_2} (M_2(\Phi^2),M_2(\tilde\Phi^2))_{[L_0,L]}
	\end{aligned}
\end{equation*}
and $D(B^{1/2})=H_d\subset H_v$.
We define also operators $K_i$ by the formulas
$$K_1(\Theta)=\left( \begin{array}{c} l\alpha_1\xi \\ \alpha_2 \theta_x \\ \alpha_1\xi_x \\0\\0\\0  \end{array} \right),\qquad\qquad K_2(\Phi_t)=\left( \begin{array}{c} \frac{\alpha_1}{\sigma} N(\Phi_{1t})\\ \frac{\alpha_2}{\nu_1} M_1(\Phi_{1t})   \end{array} \right)$$
and $C:D(C)\subset H_\theta\arr H_\theta$ 
\begin{equation*}
	C=\left(\begin{array}{cc} -\mu \partial_{xx} &    0 \\ 0 &  -\lambda \partial_{xx}    \end{array}\right)
\end{equation*}
with the domain $D(C)=H^2(0,L_0)\times H_0^1(0,L_0)$.

Next we consider the phase space 
$$H=H_d\times H_v\times H_\theta$$
and introduce an operator $A:D(A)\subset H\arr H$ of the form
\begin{equation*}
	A=\left(\begin{array}{ccc} 0 & -I  &   0\\ R^{-1}B &  0 &  R^{-1} K_1 \\     0 &  S^{-1}K_2 &  S^{-1}C   \end{array}\right)
\end{equation*}
with the domain
\begin{equation*}
	D(A)=\left\{U=(\Phi,\tilde \Phi, \Theta)\in H:  \Phi\in D(B), \tilde \Phi\in H_d , \Theta\in D(C)  \right.\}.
\end{equation*}
Thus, we can rewrite problem \eqref{Eq1}-\eqref{IC} in an abstract form
\begin{equation} \label{AEq}
	U_{t}+AU+F(U)=P(x,t),
\end{equation}
where $P(x,t)=(0,\EuScript P, G)(x,t)$ and
$F(U)=(0, F_1(U), F_2(U))$
$$
F_1(U)=F_3(U)+\left(\begin{array}{c}-\frac{\sigma l}{2}\psi^2\\  \sigma\psi(\omega_x-l\vph)+\frac{\sigma}{2}\psi^3\\
	-\sigma\psi\psi_x\\
	-\frac{\sigma l}{2}v^2\\
	\sigma v(w_x-lu)+\frac{\sigma}{2}v^3\\
	-\sigma vv_x	
\end{array}\right),\quad F_3(U)=\left(\begin{array}{c}0\\  -\alpha_1\psi\xi\\
	0\\
	0\\
	0\\
	0
\end{array}\right)
$$
$$
F_2(U)=\left(\begin{array}{c}
	\alpha_1\psi_t\psi\\0
\end{array}\right).
$$
The initial conditions have the form
\begin{equation}
	U(x,0)=U_0(x), \label{AIC}
\end{equation}
where $U_0(x)=(\Phi_0, \Phi_1,\Theta_0)$.
Using the notation  $J(\nu, \zeta, \eta)=\si(\eta_x-l\nu+\frac{\zeta^2}{2})$ we also introduce a form
\begin{equation*} \label{AForm1}
	\begin{aligned}
		\EuScript G(\Phi, \tilde \Phi)	=
		& \frac{1}{k_1} (Q_1(\Phi^1),Q_1(\tilde\Phi^1))_{[0,L_0]} + \frac{1}{\si} (J(\Phi^1),\tilde\omega_{1x}-l\tilde\varphi_1+ \psi_1\tilde \psi_1))_{[0,L_0]} + \frac{1}{\la_1} (M_1(\Phi^1),M_1(\tilde\Phi^1))_{[0,L_0]} + \\
		& \frac{1}{k_2} (Q_2(\Phi^2),Q_2(\tilde\Phi^2))_{[L_0,L]} + \frac{1}{\si} (J(\Phi^2),\tilde\omega_{2x}-l\tilde\varphi_2+ \psi_2\tilde \psi_2))_{[L_0,L]} + \frac{1}{\la_2} (M_2(\Phi^2),M_2(\tilde\Phi^2))_{[L_0,L]}.
	\end{aligned}
\end{equation*}
Note that
\begin{equation*} \label{AForm1}
	\begin{aligned}
		\EuScript G(\Phi,  \Phi)	=
		& \frac{1}{k_1} \|Q_1(\Phi^1)\|_{[0,L_0]}^2 + \frac{1}{\si} \|J(\Phi^1)\|_{[0,L_0]}^2 + \frac{1}{\la_1} \|M_1(\Phi^1)\|_{[0,L_0]}^2 + \\
		& \frac{1}{k_2} \|Q_2(\Phi^2)\|_{[L_0,L]}^2 + \frac{1}{\si} \|J(\Phi^2)\|_{[L_0,L]}^2 + \frac{1}{\la_2} \|M_2(\Phi^2)\|_{[L_0,L]}^2.
	\end{aligned}
\end{equation*}
It is easy to infer from the embedding $H^1(\Omega)\subset L_4(\Omega)$ for  $\Omega\in \mathbb R$ that there exist constants $C_i>0$, $i=1,2,3$ such that
\begin{align}
	&\|J(\nu, \zeta, \eta)\|^2\le C\|\eta_x-l\nu+\frac{\zeta^2}{2}\|^2\le  C_1(1+\|(\nu, \zeta, \eta)\|_{H^1}^4),\label{es1}\\
	&\|J(\nu, \zeta, \eta)\|^2\ge C_2\|\eta_x-l\nu\|^2-C_3\|\zeta_x\|^4.\label{es2}
\end{align}
To give the definition of a variational solution 
we introduce the space
\begin{equation*}
	F_T^0=\{V=(\Psi, \tau)\in L^2(0,T;H_d\times \tilde H_\theta):V_t=(\Psi_t, \tau_t)\in L^2(0,T;H_v\times  H_\theta),\; V(T)=0\}
\end{equation*}
of test functions.
\section{Well-posedness}
In this section we prove the existence of weak and approximate solutions to \eqref{AEq}-\eqref{AIC}.
\begin{definition}
	$U\in C(0,T;H)$ such that $U(x,0)=U_0(x)$,  is said to be  a strong solution to \eqref{AEq}-\eqref{AIC} if
	\begin{itemize}
		\item $U(t)$ lies in $D(A)$ for almost all $t$;
		\item  $U(t)$ is an absolutely continuous function with values in $H$  and $U_t \in L_1(a,b;H)$ for $0<a<b<T$;
		\item equation \eqref{AEq} is satisfied for almost all $t$.
	\end{itemize}
\end{definition}
\begin{definition}
	$U\in C(0,T;H)$ such that $U(x,0)=U_0(x)$,  is said to be  a generalized solution to \eqref{AEq}-\eqref{AIC} if there exists a sequence of strong solutions $U^{(n)}$ to \eqref{AEq}-\eqref{AIC} with the initial data $U_0^{(n)}$ and right hand sides $P^{(n)}(x,t)$ such that
	\begin{equation}
		\label{gen}
		\lim_{n\arr\infty} \max_{t\in[0,T]} ||U^{(n)}(\cdot,t)-U(\cdot,t)||_H =0.
	\end{equation}
\end{definition}

\begin{definition}
	$U(t)=(\Phi(t), \Phi_t(t) , \Theta)$ is said to be a variational (weak) solution to \eqref{AEq}-\eqref{AIC} if
	\begin{itemize}
		\item $\Phi\in L^\infty(0,T;H_d), \; \Phi_t\in L^\infty(0,T;H_v),\; \Theta\in L^\infty(0,T;H_\theta)\cap L^2(0,T;\tilde H_\theta)$;
		\item satisfy the following variational equality for all $V\in F_T^0$
		\begin{multline}\label{VEq}
			-\intl_0^T  (R\Phi_t,\Psi_t)(t)dt-\intl_0^T  (S\Theta,\tau_t)(t)dt + \int_0^T  \EuScript G(\Phi,  \Psi)dt +\int_0^T  (C^{1/2}\Theta, C^{1/2}\tau)(t)dt\\
			+\int_0^T (K_1(\Theta), \Psi)(t)dt+\int_0^T (K_2(\Phi_t), \tau)(t)dt	+\int_0^T (F_2(U), \tau)(t)dt\\+\int_0^T (F_3(U), \Psi)(t)dt= \int_0^T (P, V)(t)dt+(R\Phi_1, \Psi(0))+(S\Theta_0, \tau(0))
		\end{multline}
		\item $\Phi(x,0)=\Phi_0(x)$.
	\end{itemize}
\end{definition}
\begin{theorem}[Well-posedness] \label{th:WeakWP}
	Let
	\begin{equation}
		P\in L^2(0,T;H_v\times  H_\theta).  \label{RSmooth}
	\end{equation}
	Then for every initial data $U_0\in H$ and time interval $[0,T]$ there exists a unique generalized  solution $U(t)=(\Phi(t),\Phi_t(t),\Theta(t))\in H$ to \eqref{Eq1}-\eqref{IC}  with the following properties:
	\begin{itemize}
		\item 
		\begin{equation}
			\label{1d}
			\Phi\in C(0,T;H_d), \; \Phi_t\in C(0,T;H_v), \; \Theta\in C(0,T;H_\theta).
		\end{equation}
		\item Energy equality 
		\begin{equation}\label{EE}
			\cE(T)+\int_0^T  \|C^{1/2}\Theta\|^2dt = \cE(0)+ \int_0^T (\EuScript P(t),\Phi_t(t))dt+\int_0^T  (G(t),\Theta(t))dt
		\end{equation}
		holds, where
		\begin{equation}
			\label{tnfunc}
			\cE(t)=\cE(\Phi(t),\Phi_t(t),\Theta(t))=\cE(U(t))=\frac 12 \left[||R^{1/2}\Phi_t(t)||^2+||S^{1/2}\Theta(t)||^2 +   \EuScript G(\Phi,\Phi)\right].
		\end{equation}
		\item The solution depends continuously on initial data, i.e. if $U_n\to U_0$ in $H$, then $U(t,U_n)\to U(t;U_0)$ in $H$ for each $t$.
		\item Generalized solutions are also weak and   
		\begin{equation}
			\label{2d}
			\Phi_{tt}\in L_2(0,T;(H_v)^*)\subset L_2(0,T;[H^{-1}(0,L_0)]^3\times [H^{-1}(L_0,L)]^3).
		\end{equation}
		\item If, additionally, $U_0\in D(A)$ and
		\begin{equation}
			\pd_t P(x,t)\in L_2(0,T;H_v\times H_\theta)   \label{RAddSmooth}	
		\end{equation}
		then the generalized solution is also strong and satisfies the energy equality. 
	\end{itemize}
\end{theorem}

\begin{proof}
	The proof essentially is rather standard (see, e.g., \cite{ChuEllLa2002}), so in some parts we give only references to corresponding arguments.\\
	
	{\it Step 1. Existence and uniqueness of a local solution.} Here we use Theorem~7.2 from \cite{ChuEllLa2002}.  For the reader's convenience we formulate it below.
	\begin{theorem}[\cite{ChuEllLa2002}]\label{th:AbstrEx}
		Consider the initial value problem
		\begin{equation}\label{AbstrIVP}
			U_t+A U +F(U) = P, \quad U(0)=U_0\in H.
		\end{equation}
		Suppose that $A:D(A)\subset H \arr H$ is a maximal monotone mapping,  $0\in A0$ and $F:H\arr H$ is locally Lipschitz, i.e. there exists $L(K)>0$ such that
		\begin{equation*}
			||F(U)-F(V)||_H \le L(K)||U-V||_H, \quad ||U||_H, ||V||_H \le K.
		\end{equation*}
		If $U_0\in D(A)$, $P\in W_1^1(0,t;H)$ for all $t>0$, then there exists $t_{max}\le \infty$ such that \eqref{AbstrIVP} has a unique strong solution $U$ on $(0,t_{max})$.\\
		If $U_0\in \overline{D(A)}$, $f\in L^1(0,t;H)$ for all $t>0$, then there exists $t_{max}\le \infty$ such that \eqref{AbstrIVP} has a unique generaized solution $U$ on $(0,t_{max})$.\\
		In both cases
		\begin{equation*}
			\lim_{t\arr t_{max}}||U(t)||_H=\infty \quad \mbox{provided} \quad t_{max}<\infty.
		\end{equation*}
	\end{theorem}
	First, we need to check that $A$ is a maximal monotone operator. Monotonicity is obvious.
	To prove that $A$ is maximal as an operator from $H$ to $H$, we use Theorem 1.2 from \cite[Ch. 2]{Bar1976}. Thus, we need to prove that $Range(I+A)=H$, with $I$ being the duality map from $H$ to $H$. Let $f=(f_1,f_2,f_3)\in H$. We need to find $U=(\Phi,\tilde\Psi,\Theta)\in D(A)$ such that
	\begin{align}
		& -\tilde\Phi + \Phi=f_1\in H_d,\label{m1}\\
		& B\Phi+K_1\Theta + R\tilde\Phi = Rf_2\in H_v,\label{m2}\\
		&  K_2\tilde \Phi+C\Theta+S\Theta=Sf_3\in H_\theta.\label{m3}
	\end{align}
	Substituting $\tilde\Phi$ from \eqref{m1} to \eqref{m2},\eqref{m3} we come to
	\begin{align}
		& B\Phi+\Phi+K_1\Theta =f_1+ Rf_2\in H_v,\label{m22}\\
		& C\Theta+S\Theta+K_2\Phi=Sf_3+K_2f_1\in H_\theta.\label{m32}
	\end{align}
	By the Lax-Milgram theorem there exists a unique solution $(\Phi,\Theta)\in H_d\times D(C^{1/2})$ to \eqref{m22}, \eqref{m32} and the structure of \eqref{m1}, \eqref{m22}, \eqref{m32} implies $(\Phi,\tilde \Phi,\Theta)\in D(B)\times H_d\times D(C)=D(A)$.\\ 
	Further, we need to prove that $F$ is locally Lipschitz on $H$. The embedding $H^{1/2+\ep}(0,L)\subset C(0,L)$  implies
	\begin{equation}
		||F_j(\widetilde U)-F_j(\widehat U)||\le p(||\widetilde U||_H, ||\widehat U||_H)) ||\widetilde U-\widehat U||_{H}  \label{FLip}
	\end{equation}
	for  $j=1,2$, where $p$ is a polynomial.
	Thus, all the assumptions of Theorem \ref{th:AbstrEx} are satisfied and existence of a local strong/generalized solution is proved.\\
	{\it Step 3. Energy inequality and  global solutions.}
	It can be verified by direct calculations, that strong solutions satisfy the energy equality. Using the same arguments, as in the proof of Proposition 1.3 \cite{ChuLa2007},  we can pass to the limit and prove \eqref{EE} for generalized solutions.  \\
	Let us assume that a local generalized solution exists on a maximal interval $(0, t_{max})$, $t_{max}<\infty$. Then \eqref{EE} implies $\cE(t_{max})\le C(T,\|P\|,\cE(0))$. Therefore, due to \eqref{es1}, \eqref{es2}
	we have $||U(t_{max})||_H\le C(T,\|P\|,\|U_0\|_H)$. Thus, we arrive to a contradiction which implies $t_{max}=\infty$.
	
	{\it Step 4. Generalized solution is variational (weak).} 
	
	Evidently, \eqref{VEq} is valid for strong solutions. We can find a sequence of strong solutions $U_n$, which converges to a generalized solution $U$ strongly in $C(0,T; H)$. Using this fact we can easily pass to the  limit  in linear  terms in \eqref{VEq}. Using the embedding $H^1(0,L)\subset C(0,L)$ we get
	\begin{multline*}
		\left| \int_{0}^T(J(\Phi_n^1),\eta_{x}^1-l\nu^1+\psi_n \zeta^1)_{[0,L_0]}-(J_1(\Phi^1),\eta_{x}^1-l\nu^1+\psi \zeta^1)_{[0,L_0]}dt\right|\\\le C
		(\int_{0}^T(1+\|\Phi^1\|_{H^1}^2)\|\psi_n-\psi\|\|\zeta_{x}^1\|dt+
		\int_{0}^T\|J_1(\Phi_n^1)-J_1(\Phi^1)\|\|\eta_{x}^1-l\nu^1+\psi_n \zeta^1\|dt)\\\le C
		(
		\int_{0}^T(\|\omega_{nx}-\omega_{x}\|+\|\varphi_{n}-\varphi\|+\|\psi_n-\psi\|(\|\psi_n\|+\|\psi\|))(\|\eta_{x}^1\|+\|\nu^1\|+\|\psi_n \|\|\zeta_{x}^1\|)dt\\+\int_{0}^T(1+\|\Phi^1\|_{H^1})\|\psi_n-\psi\|\|\zeta_{x}^1\|dt)\le
		C(T,\|U_0\|_H, \|P\|)\int_{0}^T\|U_n-U\|_H \|\Psi\|_{H_d}dt\to 0,\quad n\to \infty.
	\end{multline*} 
	Arguing analogously for the second nonlinear term, we obtain the statement of the theorem. For strong solutions, one can obtain the following estimate from  equation \eqref{AbstrIVP}
	\begin{multline*}
		\int\limits_0^T\|\Phi_{tt}\|_{(H_v)^*}^2 dt\le C	\int\limits_0^T(\|B^{1/2}\Phi\|^2+\|B^{-1/2}K_1\Theta\|^2 +\|F_1(U)\|^2)dt \\\le   C(1+(\max\limits_{[0,T]}\|\Phi\|_{H_v})^4+(\max\limits_{[0,T]}\|\Theta\|_{H_\theta})^4).  
	\end{multline*}
	Due to \eqref{gen} one can infer that there exists a subsequence of strong solutions $U_n$ converging to a generalized solution $U$ such that  
	$\Phi_{ntt}$ converges weakly in $L_2(0,T;(H_v)^*)$ to a function $\bar \Phi$. Combining this with \eqref{gen} we obtain \eqref{2d}.
\end{proof}
\section{Existence of attractors.}
In this section we study long time behaviour of solutions to problem \eqref{AEq}-\eqref{AIC} in the framework of dynamical systems theory. From Theorem \ref{th:WeakWP} we have
\begin{corollary}\label{DSGen}
	Let, additionally to conditions of Theorem \ref{th:WeakWP}, $P(x,t)=P(x)$. Then \eqref{AEq}-\eqref{AIC} generates a dynamical system $(H, S_t)$ by the formula
	\begin{equation*}
		S_t(\Phi_0,\Phi_1,\Theta_0)=(\Phi(t),\Phi_t(t),\Theta(t)),
	\end{equation*}
	where $(\Phi(t),\Phi_t(t),\Theta(t))$ is the weak solution to \eqref{AEq}-\eqref{AIC} with initial data $(\Phi_0,\Phi_1,\Theta_0)$.
\end{corollary}	
To establish the existence of the attractor for this dynamical system we use Theorem \ref{abs} below, thus we need to prove the gradient property and  the asymptotic smoothness as well as the boundedness of the set of stationary points.
\subsection{Gradient structure.}
In this subsection we prove that the dynamical system generated by \eqref{AEq}-\eqref{AIC} possesses specific structure, namely, is gradient under some additional conditions on the nonlinearities.
\begin{definition}[\cite{Chueshov,CFR,CL}]\label{de:grad}
	Let $Y\subseteq X$ be a positively invariant set of  $(X,S_t)$.
	\begin{itemize}
		\item a continuous functional  $L(y)$ defined on  $Y$ is said to be a \emph{Lyapunov function} of the dynamical system  $(X,S_t)$ on the set $Y$, if
		a function  $t\mapsto L(S_ty)$ is non-increasing for any $y\in Y$.
		\item
		the Lyapunov function  $L(y)$ is said to be \emph{strict} on $Y$,
		if the equality  $L(S_{t}y)=L(y)$  \emph{for all} $t>0$ implies  $S_{t}y=y$ for all $t>0$;
		\item
		A dynamical system $(X,S_t)$ is said to be \emph{gradient}, if it possesses a strict Lyapunov function on the whole phase space $X$.
	\end{itemize}
\end{definition}
The following result holds true.
\begin{theorem}\label{th:grad}
	Let, assumptions of Corollary \ref{DSGen} hold
	and
	\begin{equation}
		\label{rightside}
		G(x)=(h_1(x), h_2(x))=(0,0).
	\end{equation}
	Then the dynamical system $(H, S_t)$ is gradient.
\end{theorem}
\begin{proof}
	We use as a Lyapunov function
	\begin{equation}\label{lap}
		L(\Phi,\Theta)=
		\cE(t)-(\EuScript P,\Phi(t)),
	\end{equation}
	where $\cE(t)$ is defined in \eqref{tnfunc}.
	Energy equality \eqref{EE} implies that $L(t)$ is non-increasing. The equality $L(t)=L(0)$  yields  $\Theta(t)\equiv 0$ on $[0,T]$.  We need to prove that  $\Phi(t)\equiv const$, which is equivalent to $\Phi(t+h)-\Phi(t)=0$ for every $h>0$. In what follows we denote $\Phi(t+h)-\Phi(t)=\oPhi(t)=(\ovph, \overline{\psi},\oom,\ou,\ov,\ow)(t)$.
	Analogously, substituting $V=(0,0,0,0,0,0,0,\tau_2)\in C_0^\infty((0,L_0)\times (0,T))$  into \eqref{VEq}  as a test function, we obtain  $\psi_{tx}=0$ in the sense of distributions.	Similar to regular functions,  if the partial derivative of a distribution equals to zero, then the distribution does not depend on the corresponding variable (see \cite[Ch. 7]{Kan2004}, Example 2). That is,
	\begin{equation*}
		{\psi}_x=c_1(x)\times 1(t) \qquad \mbox{ in the sense of  distributions} .
	\end{equation*}
	However, Theorem \ref{th:WeakWP} implies that  
	${\psi}_x$ is a regular distribution, thus, we can treat the equality above as the equality almost everywhere and obtain that ${\psi_x}$ does not depend on $t$. Taking into account boundary conditions, we get that ${\psi}$ also does not depend on $t$, which means $\overline{\psi}=0$.\\ 
	If we subsitute  a test function $V=((\Psi^1, 0),0)\in C_0^\infty((0,L_0)\times (0,T))$ into \eqref{VEq} we obtain
	\begin{equation}\label{VEqdif}
		\intl_0^T  (R\Phi_{tt}^1,\Psi^1)(t)dt + \int_0^T  \EuScript G(\Phi^1,  \Psi^1)dt=	\intl_0^T  (\EuScript P(x),\Psi(t))dt.
	\end{equation}
	Then \eqref{VEqdif} implies
	\begin{align} 
		& \ovph_x+l\oom=0 \qquad&\mbox{ almost everywhere on } (0,L_0)\times (0,T),  \label{GrSpace} \\
		& \rho_1\oom_{tt} - l\si_1(\oom_x-l\ovph)_x =0 \qquad&\mbox{ in the sense of distributions}, \label{GrOmega} \\
		& \rho_1\ovph_{tt} - \si_1(\oom_x-l\ovph)=0,  \qquad&\mbox{ in the sense of distributions  }. \label{GrPhi}
	\end{align}
	These equalities imply
	\begin{equation}\label{GrTime}
		\oom_{tt}=0\quad \mbox{ in the sense of distributions}.
	\end{equation}
	Consequently, $\oom$ does not depend on $t$, moreover, since $\oom_t\in L_2(0,T)$ for  $x\in (0,L_0)$,
	\begin{equation*}
		\oom(x,t)=\oom(x,0)+\int_0^t \oom_\tau (x)d\tau = \oom(x,0) +t\oom_t(x).
	\end{equation*}
	Since $||\oom(\cdot,t)||\le C$ for all $t\in\R_+$, $\oom_t(x)$ must be zero. Thus,
	\begin{equation*}\label{GrOmegaConst}
		\oom(x,t)=c_2(x),
	\end{equation*}
	which together with \eqref{GrSpace} implies
	\begin{align*}
		& \ovph_x=-lc_2(x), \\
		& \ovph(x,t)= \ovph(0,t) - l\intl_0^x  c_2(y)dy = c_3(x), \\
		& \ovph_{tt} =0.
	\end{align*}
	The last equality together with \eqref{GrSpace}, \eqref{GrPhi} and the boundary conditions  give us that $\ovph, \oom$ are solutions to the following Cauchy problem (with respect to $x$):
	\begin{align*}
		&\oom_x = l\ovph,\\
		&\ovph_x = -l\oom, \\
		&\oom(0,t)=\ovph(0,t)=0.
	\end{align*}
	Consequently, $\oom\equiv\ovph\equiv 0$.\\
	{\it Step 2.} Let us prove, that $\bar u\equiv \bar v\equiv \bar w\equiv 0$. It is easy to see that  $(\ou, \ov, \ow)$ is a weak solution to a linear overdetermined problem
	on $(0,T)\times (L_0,L)$
	\begin{align}
		& \rho_2\ou_{tt}-k_2\ou_{xx} +g_u(\pd_x \oPhi^2, \oPhi^2,x,t) =0, \label{WE1} \\
		& \beta_2\ov_{tt} -\nu_2 \ov_{xx} + g_v(\pd_x \oPhi^2, \oPhi^2,x,t)=0,\\
		& \rho_2\ow_{tt}- \si \ow_{xx} + g_w(\pd_x \oPhi^2, \oPhi^2,x,t) =0, \\
		& \ou(L_0,t)=\ov(L_0,t)=\ow(L_0,t)=0,\\
		& \ou(L,t)=\ov(L,t)=\ow(L,t)=0, \\
		& {\ou_x(L_0,t)=0}, \quad
		{\ov_x(L_0,t)=0, \quad \ow_x(L_0,t)=0}, \label{WEBC2}\\
		& \oPhi^2(x,0)=\Phi^2(x,h)-\Phi^2_0, \quad \oPhi^2_t(x,0)=\Phi^2_t(x,h)-\Phi^2_1, \label{WEIC}
	\end{align}
	where $g_u, g_v, g_w$ are linear combinations of $u_x,v_x,w_x, u,v,w$ with the coefficients depending on $x$ and $t$.  The $L^2$-regularity of $u_x,v_x, w_x$ on the boundary for solutions to a linear wave equation  was established in \cite{LaTri1983}, thus, boundary conditions \eqref{WEBC2} make sense.\\
	It is easy to generalize the Carleman estimate \cite[Th. 8.1]{TriYao2002} for the case of the system of the wave equations.
	\begin{theorem}[\cite{TriYao2002} ]
		For the solution to problem \eqref{WE1}-\eqref{WEIC} the following estimate holds:
		\begin{equation*}
			\int_0^T  [|\ou_x|^2+|\ov_x|^2+|\ow_x|^2](L_0,t) dt\ge C(\hat E(0)+\hat E(T)),
		\end{equation*}
		where
		\begin{equation*}
			\hat E(t)=\frac 12 \left(||\ou_t(t)||^2+ ||\ov_t(t)||^2 +||\ow_t(t)||^2 + ||\ou_x(t)||^2 +||\ov_x(t)||^2+||\ow_x(t)||^2 \right).
		\end{equation*}
	\end{theorem}
	Therefore, if conditions  \eqref{WEBC2} hold true, then $\ou=\ov=\ow=0$.
	The theorem is proved.
\end{proof}
\subsection{Asymptotic smoothness.}
\begin{definition}[\cite{Chueshov,CFR,CL}]
	A dynamical system $(X,S_t)$ is said to be asymptotically smooth
	if for any  closed bounded set $B\subset X$ that is positively invariant ($S_tB\subseteq B$)
	one can find a compact set $\cK=\cK(B)$ which uniformly attracts $B$, i.~e.
	$\sup\{{\rm dist}_X(S_ty,\cK):\ y\in B\}\to 0$ as $t\to\infty$.
\end{definition}
In order to prove  the asymptotical smoothness of the system
considered we rely on the  compactness criterion due to
\cite{Khanmamedov}, which is recalled below in an abstract version
formulated in \cite{CL}.
\begin{theorem}{\cite{CL}} \label{theoremCL} Let $(S_t, H)$ be a dynamical system on a complete metric
	space $H$ endowed with a metric $d$. Assume that for any bounded positively invariant
	set $B$ in $H$ and for any $\varepsilon>0$ there exists $T = T (\varepsilon, B)$ such that
	\begin{equation}
		\label{te}
		d(S_T y_1, S_T y_2) \le \varepsilon+ \Psi_{\varepsilon,B,T} (y_1, y_2), y_i \in B ,
	\end{equation}
	where $\Psi_{\varepsilon,B,T} (y_1, y_2)$ is a function defined on $B \times B$ such that
	\[
	\liminf\limits_{m\to\infty}\liminf\limits_{n\to\infty}\Psi_{\varepsilon,B,T} (y_n, y_m) = 0
	\]
	for every sequence ${y_n} \in B$. Then $(S_t, H)$ is an asymptotically smooth dynamical
	system.
\end{theorem}
\begin{theorem}
	\label{th:AsSmooth}
	Let the assumptions of Theorem \ref{th:WeakWP} and \eqref{rightside}
	hold. Let, moreover,
	\begin{equation}
		\rho_1>\rho_2, \beta_1>\beta_2,\; k_1\le k_2, \,\nu_1\le \nu_2, \,\frac{\rho_1}{\beta_1}=\frac{k_1}{\nu_1}.  \label{c1}
	\end{equation}
	Then the dynamical system  $(H, S_t)$ generated by problem \eqref{AEq}-\eqref{AIC} is asymptotically smooth.
\end{theorem}
\begin{proof}
	In this proof we perform all the calculations for strong solutions and then pass to the limit in the final estimate to justify it for weak solutions.
	Let us consider strong solutions $\hat U(t)=(\hat\Phi(t), \hat\Phi_t(t), \hat\Theta(t))$ and $\tilde U(t)=(\tilde\Phi(t), \tilde\Phi_t(t), \tilde\Theta(t))$  to problem \eqref{AEq}-\eqref{AIC} with initial conditions $\hat U_0=(\hat \Phi_0, \hat \Phi_1, \hat\Theta_0)$ and $\tilde U_0=(\tilde \Phi_0, \tilde \Phi_1, \tilde\Theta_0)$ lying in a ball, i.e.  there exists $R>0$ such that
	\begin{equation}
		\label{inboun}
		\|\tilde U_0\|_H+\|\hat U_0\|_H\le R.
	\end{equation}
	Denote $U(t)=\tilde U(t)-\hat U(t)$ and  $U_0=\tilde U_0-\hat U_0$. Obviously, $U(t)$ is a weak solution to the problem
	\begin{align}
		\begin{split} \rho_1\vph_{tt}-k_1(\vph_x+\psi+l\om)_x - l\si(\om_x-l\vph)-\frac{\sigma l}{2}\psi(\hat\psi+\tilde\psi)+l\alpha_1\xi=0,\;x\in (0,L_0), t>0,\label{eq1}
		\end{split}\\
		\begin{split} \be_1\psi_{tt} -\nu_1 \psi_{xx} +k_1(\vph_x+\psi+l\om) +\alpha_2\theta_x-\frac{\alpha_1}{2}\xi(\tilde\psi+\hat\psi)-\frac{\alpha_1}{2}\psi(\tilde\xi+\hat\xi)\qquad\qquad\qquad\\\qquad
			+\frac{\sigma}{2}(\tilde\psi+\hat\psi)(\omega_x-l\varphi+\frac12(\tilde\psi+\hat\psi)\psi)+\frac{\sigma}{2}\psi(\tilde\omega_x-l\tilde\varphi+\hat\omega_x-l\hat\varphi+\frac{\hat\psi^2}{2}+\frac{\tilde\psi^2}{2})=0,\label{eq2}
		\end{split}\\
		& \rho_1\om_{tt}- \si(\om_x-l\vph) _x+lk_1(\vph_x+\psi+l\om)-\frac{\sigma}{2}\psi(\tilde\psi_x+\hat\psi_x)-\frac{\sigma}{2}(\tilde\psi+\hat\psi)\psi_x+\alpha_1\xi_x=0, \label{eq3}\\
		&\gamma\xi_t-\mu\xi_{xx} +\alpha_1(\omega_x-l\vph)_t+\frac{\alpha_1}{2}(\tilde\psi_t+\hat\psi_t)\psi+\frac{\alpha_1}{2}\psi_t(\hat\psi+\tilde\psi)=0,\label{eq31}\\
		&\delta\theta_t-\lambda\theta_{xx}+\alpha_2\psi_{tx}=0\label{eq32}
	\end{align}
	and
	\begin{align}
		& \rho_2u_{tt}-k_2(u_x+v+lw)_x - l\si(w_x-lu) -\frac{\sigma    l}{2}v(\hat v+\tilde v)=0,\qquad x\in (L_0,L), t>0, \label{eq4}\\
		\begin{split}
			\be_2v_{tt}-\nu_2 v_{xx}+k_2(u_x+v+lw)+
			\frac{\sigma}{2}(\tilde v+\hat v)(w_x-lu+\frac12(\tilde v+\hat v)v)\qquad\qquad\qquad\qquad\\+\frac{\sigma}{2}v(\tilde w_x-l\tilde v+\hat w_x-l\hat v+\frac{\hat v^2}{2}+\frac{\tilde v^2}{2}) =0,\label{eq5}
		\end{split}
		\\
		& \rho_2w_{tt}- \si(w_x-lu)_x+lk_2(u_x+v+lw)-\frac{\sigma}{2}v(\tilde v_x+\hat v_x)-\frac{\sigma}{2}(\tilde v+\hat v)v_x=0, \label{eq6}
	\end{align}
	with boundary conditions \eqref{BC}--\eqref{TC4} and the initial conditions $U(0)=\tilde U_0-\hat U_0$.
	It is easy to see by the energy argument that
	\begin{equation}
		\label{En1}
		E(U(T))+ \int\limits_t^T \int\limits_0^{L_0}(\mu \xi_x^2+\lambda\theta_x^2) dx ds=E(U(t))+ \int\limits_t^T  H(\hat U(s),\tilde U(s))  ds,
	\end{equation}
	where
	\begin{multline*}
		\label{H}
		H(\hat U(t),\tilde U(t))=\frac{\sigma l}{2}\int\limits_0^{L_0}\psi(\hat \psi+\tilde\psi)\vph_t dx+\frac{\alpha_1}{2}\int\limits_0^{L_0}\psi(\tilde\xi+\hat\xi)\psi_t dx-\frac{\alpha_1}{2}\int\limits_0^{L_0}\psi(\tilde\psi_t+\hat\psi_t)\xi dx\\+	\frac{\sigma l}{2}\int\limits_0^{L_0}(\tilde\psi+\hat\psi)\psi_x\omega_t dx-	\frac{\sigma l}{2}\int\limits_0^{L_0}\left(\tilde\omega_x-l\tilde\vph+\hat\omega_x-l\hat\vph+\frac{\tilde\psi^2}{2}+\frac{\hat\psi^2}{2}\right)\psi\psi_t dx\\ +
		\frac{\sigma l}{2}\int\limits_0^{L_0}(\tilde\psi_x+\hat\psi_x)\psi\omega_t dx-	\frac{\sigma l}{2}\int\limits_0^{L_0}\left(\omega_x-l\vph+\frac{1}{2}\psi(\tilde\psi+\hat\psi)\right)(\tilde\psi+\hat\psi)\psi_t dx\\
		+
		\frac{\sigma l}{2}\int\limits_0^{L_0}(\tilde v_x+\hat v_x)v u_t dx+\frac{\sigma l}{2}\int\limits_{L_0}^Lv(\hat v+\tilde v)u_t dx
		+	\frac{\sigma l}{2}\int\limits_{L_0}^L(\tilde v+\hat v)v_x w_t dx-	\\\frac{\sigma l}{2}\int\limits_{L_0}^L\left(\tilde\omega_x-l\tilde\vph+\hat\omega_x-l\hat\vph+\frac{\tilde v^2}{2}+\frac{\hat v^2}{2}\right)v v_t dx -	\frac{\sigma l}{2}\int\limits_{L_0}^L\left(w_x-lu+\frac{1}{2}v(\tilde v+\hat v)\right)(\tilde v+\hat v)v_t dx	,
	\end{multline*}
	and
	\begin{equation*}
		\label{E}
		E(t)=E_1(t)+E_2(t),
	\end{equation*}
	where
	\begin{multline*}
		E_1(t)=\rho_1\int\limits_0^{L_0} \omega_t^2dx dt+\rho_1 \int\limits_0^{L_0}\vph_{t}^2 dx dt+\beta_1 \int\limits_0^{L_0}\psi_{t}^2 dx+\sigma\int\limits_0^{L_0} (\omega_x-l\vph)^2 dx+\\+k_1\int\limits_0^{L_0} (\vph_x+\psi+l\om)^2 dx+\nu_1\int\limits_0^{L_0} \psi_x^2 dx+\gamma \int\limits_0^{L_0} \xi^2 dx+\delta \int\limits_0^{L_0} \theta^2 dx
	\end{multline*}
	and
	\begin{multline*}
		\label{E2}
		E_2(t)=\rho_2\int\limits_0^{L_0} w_t^2dx dt+\rho_2 \int\limits_0^{L_0}u_{t}^2 dx dt+\beta_2 \int\limits_0^{L_0}v_{t}^2 dx+\sigma\int\limits_0^{L_0} (w_x-lu)^2 dx+\\+k_2\int\limits_0^{L_0} (u_x+v+lw)^2 dx+\nu_2\int\limits_0^{L_0} v_x^2 dx.
	\end{multline*}
	Integrating in \eqref{En1}  over the interval $(0,T)$ we come to
	\begin{equation}
		\label{En2}
		TE(U(T))+\int\limits_0^T  \int\limits_t^T \int\limits_0^{L_0} (\mu \xi_x^2+\lambda\theta_x^2) dx ds dt\\=\int\limits_0^T  E(U(t)) dt+ \int\limits_0^T \int\limits_t^T  H(\hat U(s),\tilde U(s))  ds dt.
	\end{equation}
	Now we  estimate the first term in the right-hand side of \eqref{En2}. In what follows we present formal estimates which can be performed on strong solutions.\\
	{\it Step 1.} 
	We introduce the following notation 
	\begin{multline}
		\label{lot}
		lot=\max\limits_{t\in[0,T]} (\| \vph(\cdot,t)\|_{H^{1-\epsilon}}^2+\| \psi(\cdot,t)\|_{H^{1-\epsilon}}^2+ \|\omega(\cdot,t)\|_{H^{1-\epsilon}}^2+\| \psi_t(\cdot,t)\|_{H^{-\epsilon}}^2+\| \omega_t(\cdot,t)\|_{H^{-\epsilon}}^2\\
		+\| u(\cdot,t)\|_{H^{1-\epsilon}}^2+ \|v(\cdot,t)\|_{H^{1-\epsilon}}^2+\|w(\cdot,t)\|_{H^{1-\epsilon}}^2),\quad 	0<\epsilon<1/2.
	\end{multline}
	It follows from \eqref{1d}, \eqref{2d} and the Aubin's compactness theorem \cite{Aubin} that \eqref{lot} is a compact seminorm in $H$.
	It follows from Theorem \ref{th:WeakWP} that if \eqref{inboun} holds, then 
	\begin{multline}
		\label{bounsolu}
		\max\limits_{[0,T]}\left(\int\limits_0^{L_0}\tilde \psi_x^2(x,t)dx +\int\limits_0^{L_0}\hat \psi_x^2(x,t)dx+\int\limits_{L_0}^L\tilde v_x^2(x,t)dx +\int\limits_{L_0}^L\hat v_x^2(x,t)dx+
		\int\limits_0^{L_0}(\tilde \omega_x-l\tilde\vph)^2(x,t)dx\right.\\\left. +\int\limits_0^{L_0}(\hat \omega_x-l\hat\vph)^2(x,t)dx+\int\limits_{L_0}^L(\tilde w_x-l\tilde u)^2(x,t)dx +\int\limits_{L_0}^L(\hat w_x-l\hat u)^2(x,t)dx\right)\le C(R).
	\end{multline}
	Let us note that integrating by parts we get
	$$
	\int\limits_{0}^{L_0}\psi_x\psi dx=\frac12 \psi^2(L_0),
	$$
	therefore, using that $(H^\epsilon)^*=(H_0^\epsilon)^* =H^{-\epsilon}$ for $0<\epsilon<1/2$ we get
	$$
	\psi^2(L_0)\le c\|\psi\|_{H^{1-\epsilon}} \|\psi\|_{H^{\epsilon}}\le c \|\psi\|_{H^{1-\epsilon}}^2.
	$$
	Arguing in the same way for other unknown functions in our problem we infer that
	\begin{equation}
		\label{lotest}
		\int_0^T(\psi^2(L_0)+\vph^2(L_0)+\omega^2(L_0))dt\le C(T)lot.
	\end{equation}
	Next we multiply equation \eqref{eq1} by $\vph_x$, equation \eqref{eq2} by $\psi_x$, equation \eqref{eq3} by $\omega_x$,  then perform integration by parts with respect to both $t$ and $x$. Taking into account that for any function $f\in H^2(0,L_0)$
	$$\int\limits_{0}^{L_0}f_{xx}f_x dx=\frac12 (f_x^2(L_0)-f_x^2(0))$$ 
	we obtain 
	\begin{multline}
		\label{l01}
		\min\{\rho_1, \beta_1\}\int\limits_0^T(\vph_t^2(L_0)+\psi_t^2(L_0)+\omega_t^2(L_0) )dt\\+ \frac{1}{2}\min\{\nu_1, \sigma, k_1\}\int\limits_0^T((\vph_x+\psi+l\omega)^2(L_0)+\psi_{x}^2(L_0)+ (\om_x-l\vph)^2)(L_0) )dt\\\le \max\{\nu_1, \sigma, k_1\}\int\limits_0^T((\vph_x+\psi+l\omega)^2(0)+\psi_{x}^2(0)+ (\om_x-l\vph)^2)(0))dt+C(E(0)+E(T))+C(R,T)lot \\
		C\int\limits_0^T\int\limits_{0}^{L_0}(\psi_{x}^2+ (\om_x-l\vph)^2 +(\vph_x+\psi+l\omega)^2+\xi_x^2+\theta_x^2)dx dt.
	\end{multline}
	{\it Step 2.} 
	Let us introduce a function $g(x)$  defined on $[0,L]$ with the following properties:
	\begin{enumerate}
		\item[(g1)] $g\in C^1(0,L)$ and   $g_1=\max\limits_{[0, L]}|g'|$,
		\item[(g2)] $g'>0$ on $[\hat L, L]$, where $0<\hat L<L_0$ and $g_0=\min\limits_{[\hat L, L]}g'>0$,
		\item[(g3)] $g(0)>0$, $g(L)<0$.
	\end{enumerate}
	We multiply equation \eqref{eq1} by $g\vph_x$, \eqref{eq2} by $g\psi_x$, and \eqref{eq3} by $g\omega_x$  and sum up the results and integrate over $[0,L_0]$ with respect to $x$ and integrate by parts over $[0,T]$ with respect to $t$
	\begin{multline*}
		-\rho_1 \int\limits_0^T\int\limits_0^{L_0}\vph_tg\vph_{tx} dx dt-\beta_1 \int\limits_0^T\int\limits_0^{L_0}\psi_tg\psi_{tx} dx dt -\rho_1 \int\limits_0^T\int\limits_0^{L_0} \omega_tg\omega_{tx}   dx dt\\-
		k_1 \int\limits_0^T \int\limits_0^{L_0}(\vph_x+\psi+l\omega)_xg\vph_x  dx dt
		-l\sigma \int\limits_0^T \int\limits_0^{L_0}(\omega_x-l\vph)g\vph_x   dx dt-\frac{\sigma l}{2} \int\limits_0^T \int\limits_0^{L_0} \psi(\tilde\psi+\hat\psi)g\vph_x   dx dt\\
		+\alpha_1 l \int\limits_0^T \int\limits_0^{L_0} \xi g\vph_x   dx dt
		- \nu_1\int\limits_0^T \int\limits_0^{L_0} \psi_{xx}g\psi_x dx dt+k_1 \int\limits_0^T \int\limits_0^{L_0}(\vph_x+\psi+l\omega)g\psi_x  dx dt\\
		+\alpha_2  \int\limits_0^T \int\limits_0^{L_0} \theta_x g\psi_x   dx dt-\frac{\alpha_1}{2}\int\limits_0^T \int\limits_0^{L_0} \xi(\hat\psi+\tilde\psi) g\psi_x   dx dt-\frac{\alpha_1}{2}\int\limits_0^T \int\limits_0^{L_0}\psi(\hat\xi+\tilde\xi) g\psi_x   dx dt\\+\frac{\sigma}{2}\int\limits_0^T \int\limits_0^{L_0} (\hat\psi+\tilde\psi)\left(\omega_x-l\vph+\frac{1}{2}(\hat\psi+\tilde\psi)\psi \right) g\psi_x   dx dt+\frac{\sigma}{2}\int\limits_0^T \int\limits_0^{L_0}\psi(\hat\omega_x-l\hat\vph+\tilde\omega_x-l\tilde\vph+\frac{\hat\psi^2}{2}+\frac{\tilde\psi^2}{2}) g\psi_x   dx dt\\
		-\sigma  \int\limits_0^T\int\limits_0^{L_0}(\omega_x-l\vph)_xg\omega_x   dx dt
		+lk_1\int\limits_0^T \int\limits_0^{L_0}(\vph_x+\psi+l\omega)g\omega_x  dx dt+\alpha_1  \int\limits_0^T \int\limits_0^{L_0} \xi_x g\omega_x   dx dt
	\end{multline*}
	\begin{multline*}
		-\frac{\sigma}{2}\int\limits_0^T \int\limits_0^{L_0}\psi(\hat\psi_x+\tilde\psi_x) g\omega_x   dx dt	-\frac{\sigma}{2}\int\limits_0^T \int\limits_0^{L_0}\psi_x(\hat\psi+\tilde\psi) g\omega_x   dx dt
		\\=\rho_1  \int\limits_0^{L_0} \vph_t(x,T) g \vph_{x}(x,T)  dx+\beta_1  \int\limits_0^{L_0} \psi_t(x,T) g \psi_{x}(x,T)  dx+\rho_1 \int\limits_0^{L_0} \omega_t(x,T)g \omega_x(x,T) dx\\
		-\rho_1  \int\limits_0^{L_0} \vph_1 g \vph_{0x}  dx-\beta_1  \int\limits_0^{L_0} \psi_1 g \psi_{0x}  dx-\rho_1 \int\limits_0^{L_0} \omega_1 g \omega_{0x} dx.
	\end{multline*}
	After integration by parts with respect to $x$ and taking into account boundary conditions we arrive at
	\begin{multline*}
		\frac{\rho_1}{2} \int\limits_0^T \int\limits_0^{L_0}\vph_t^2g' dx dt+\frac{\rho_1|g(L_0)|}{2} \int\limits_0^T\vph_t^2(L_0) dt+\frac{\beta_1}{2} \int\limits_0^T\int\limits_0^{L_0}\psi_t^2g' dx dt+\frac{\beta_1|g(L_0)|}{2} \int\limits_0^T\psi_t^2(L_0) dt\\+\frac{\rho_1}{2}  \int\limits_0^T\int\limits_0^{L_0} \omega_t^2g' dx dt+\frac{\rho_1|g(L_0)|}{2} \int\limits_0^T\omega_t^2(L_0) dt+
		k_1 \int\limits_0^T (\vph_x+\psi+l\omega)(L_0)g(L_0)\psi(L_0)  dt\\+k_1l \int\limits_0^T (\vph_x+\psi+l\omega)(L_0)g(L_0)\omega(L_0)  dt+\frac{k_1|g(L_0)|}{2} \int\limits_0^T (\vph_x+\psi+l\omega)^2(L_0)  dt+\frac{k_1g(0)}{2} \int\limits_0^T (\vph_x+\psi+l\omega)^2(0)  dt\\+\frac{k_1}{2} \int\limits_0^T \int\limits_0^{L_0}(\vph_x+\psi+l\omega)^2 g'  dx dt
		+l\sigma \int\limits_0^T \int\limits_0^{L_0}(\omega_x-l\vph)g'\vph  dx dt
		-l\sigma \int\limits_0^T (\omega_x-l\vph)(L_0)g(L_0)\vph(L_0)   dt\\
		+\frac{\sigma}{2} \int\limits_0^T \int\limits_0^{L_0}(\omega_x-l\vph)^2g'  dx dt+\frac{\sigma |g(L_0)|}{2} \int\limits_0^T (\omega_x-l\vph)^2(L_0)  dt+\frac{\sigma g(0)}{2} \int\limits_0^T (\omega_x-l\vph)^2(0)  dt		\\
		+\frac{\nu_1}{2}\int\limits_0^T \int\limits_0^{L_0} \psi_x^2 g' dx dt+\frac{\nu_1|g(L_0)|}{2}\int\limits_0^T  \psi_x^2(L_0) dt+\frac{\nu_1g(0)}{2}\int\limits_0^T  \psi_x^2(0) dt
		\\
		+\alpha_1 l \int\limits_0^T \int\limits_0^{L_0} \xi g\vph_x   dx dt
		+\alpha_2  \int\limits_0^T \int\limits_0^{L_0} \theta_x g\psi_x   dx dt-\frac{\alpha_1}{2}\int\limits_0^T \int\limits_0^{L_0} \xi(\hat\psi+\tilde\psi) g\psi_x   dx dt-\frac{\alpha_1}{2}\int\limits_0^T \int\limits_0^{L_0}\psi(\hat\xi+\tilde\xi) g\psi_x   dx dt\\+\frac{\sigma}{2}\int\limits_0^T \int\limits_0^{L_0} (\hat\psi+\tilde\psi)\left(-l\vph+\frac{1}{2}(\hat\psi+\tilde\psi)\psi \right) g\psi_x   dx dt+\frac{\sigma}{2}\int\limits_0^T \int\limits_0^{L_0}\psi(\hat\omega_x-l\hat\vph+\tilde\omega_x-l\tilde\vph+\frac{\hat\psi^2}{2}+\frac{\tilde\psi^2}{2}) g\psi_x   dx dt
	\end{multline*}
	\begin{multline}
		\label{gest1}
		+\alpha_1  \int\limits_0^T \int\limits_0^{L_0} \xi_x g\omega_x   dx dt-\frac{\sigma l}{2} \int\limits_0^T \int\limits_0^{L_0} \psi(\tilde\psi+\hat\psi)g\vph_x   dx dt
		-\frac{\sigma}{2}\int\limits_0^T \int\limits_0^{L_0}\psi(\hat\psi_x+\tilde\psi_x) g\omega_x   dx dt
		\\=\rho_1  \int\limits_0^{L_0} \vph_t(x,T) g \vph_{x}(x,T)  dx+\beta_1  \int\limits_0^{L_0} \psi_t(x,T) g \psi_{x}(x,T)  dx+\rho_1 \int\limits_0^{L_0} \omega_t(x,T)g \omega_x(x,T) dx\\
		-\rho_1  \int\limits_0^{L_0} \vph_1 g \vph_{0x}  dx-\beta_1  \int\limits_0^{L_0} \psi_1 g \psi_{0x}  dx-\rho_1 \int\limits_0^{L_0} \omega_1 g \omega_{0x} dx.
	\end{multline}
	If we now multiply equation \eqref{eq4} by $g u_x$, \eqref{eq5} by $g v_x$, and  \eqref{eq6} by $g w_x$  and sum up the results and integrate over $[0,L_0]$ with respect to $x$ and integrate by parts over $[0,T]$ with respect to $t$
	\begin{multline}
		\label{gest2}
		\frac{\rho_2}{2} \int\limits_0^T \int\limits_{L_0}^L u_t^2g' dx dt-\frac{\rho_2|g(L_0)|}{2} \int\limits_0^T u_t^2 dt+\frac{\beta_2}{2} \int\limits_0^T\int\limits_{L_0}^L v_t^2g' dx dt-\frac{\beta_2|g(L_0)|}{2} \int\limits_0^T v_t^2(L_0) dt\\+\frac{\rho_2}{2}  \int\limits_0^T\int\limits_{L_0}^Lw_t^2g' dx dt-\frac{\rho_2|g(L_0)|}{2} \int\limits_0^T w_t^2(L_0) dt-
		k_2 \int\limits_0^T (u_x+v+lw)(L_0)g(L_0)v(L_0)  dt\\-k_2l \int\limits_0^T (u_x+v+lw)(L_0)g(L_0)w(L_0)  dt-\frac{k_2|g(L_0)|}{2} \int\limits_0^T (u_x+v+lw)^2(L_0)  dt+\frac{k_2|g(L)|}{2} \int\limits_0^T (u_x+v+lw)^2(L)  dt\\+\frac{k_2}{2} \int\limits_0^T \int\limits_{L_0}^L(u_x+v+lw)^2 g'  dx dt
		+l\sigma \int\limits_0^T \int\limits_{L_0}^L(w_x-lu)g'u  dx dt
		+l\sigma \int\limits_0^T (w_x-lu)(L_0)g(L_0)u(L_0)   dt\\
		+\frac{\sigma}{2} \int\limits_0^T \int\limits_{L_0}^L(w_x-lu)^2g'  dx dt-\frac{\sigma |g(L_0)|}{2} \int\limits_0^T (w_x-lu)^2(L_0)  dt+\frac{\sigma |g(L)|}{2} \int\limits_0^T (w_x-lu)^2(L)  dt		\\
		+\frac{\nu_2}{2}\int\limits_0^T \int\limits_{L_0}^L v_x^2 g' dx dt-\frac{\nu_2|g(L_0)|}{2}\int\limits_0^T  v_x^2(L_0) dt+\frac{\nu_2|g(L)|}{2}\int\limits_0^T  v_x^2(L) dt
		\\
		+\frac{\sigma}{2}\int\limits_0^T \int\limits_{L_0}^L (\hat v+\tilde v)\left(-lu+\frac{1}{2}(\hat v+\tilde v)v \right) g v_x   dx dt+\frac{\sigma}{2}\int\limits_0^T \int\limits_{L_0}^L v(\hat w_x-l\hat u+\tilde w_x-l\tilde u+\frac{\hat v^2}{2}+\frac{\tilde v^2}{2}) g v_x   dx dt\\
		-\frac{\sigma l}{2} \int\limits_0^T \int\limits_{L_0}^L v(\tilde v+\hat v)g u_x   dx dt
		-\frac{\sigma}{2}\int\limits_0^T \int\limits_{L_0}^L v(\hat v_x+\tilde v_x) g w_x   dx dt
		\\=\rho_2  \int\limits_{L_0}^L u_t(x,T) g u_{x}(x,T)  dx+\beta_2  \int\limits_{L_0}^L v_t(x,T) g v_{x}(x,T)  dx+\rho_2 \int\limits_{L_0}^L w_t(x,T)g w_x(x,T) dx\\
		-\rho_2  \int\limits_{L_0}^L u_1 g u_{0x}  dx-\beta_2  \int\limits_{L_0}^L v_1 g v_{0x}  dx-\rho_2 \int\limits_{L_0}^L w_1 g w_{0x} dx.
	\end{multline}
	Summing up \eqref{gest1} and \eqref{gest2} and taking into account \eqref{l01} we obtain that the following estimate holds  for deliberately small $\hat \varepsilon>0$
	\begin{multline}
		\label{gest3}
		\frac{\rho_1 g_0}{2} \int\limits_0^T \int\limits_{\hat L}^{L_0}\vph_t^2 dx dt+\frac{\beta_1 g_0}{2} \int\limits_0^T\int\limits_{\hat L}^{L_0}\psi_t^2 dx dt+\frac{\rho_1g_0}{2}  \int\limits_0^T\int\limits_{\hat L}^{L_0} \omega_t^2 dx dt
		+\frac{\rho_1 g_0}{2} \int\limits_0^T \int\limits_{L_0}^Lu_t^2 dx dt+\frac{\beta_1 g_0}{2} \int\limits_0^T\int\limits_{L_0}^Lv_t^2 dx dt\\+\frac{\rho_1g_0}{2}  \int\limits_0^T\int\limits_{L_0}^L w_t^2 dx dt		
		+\frac{(\rho_1-\rho_2)|g(L_0)|}{2} \int\limits_0^T\vph_t^2(L_0) dt+\frac{(\beta_1-\beta_2)|g(L_0)|}{2} \int\limits_0^T\psi_t^2(L_0) dt\\+\frac{(\rho_1-\rho_2)|g(L_0)|}{2} \int\limits_0^T\omega_t^2(L_0) dt
		+	\left(\frac{k_1(k_2-k_1)|g(L_0)|}{2k_2} -\hat\varepsilon\right)\int\limits_0^T (\vph_x+\psi+l\omega)^2(L_0)  dt\\+\left(\frac{k_1g(0)}{2}-\hat\varepsilon\right) \int\limits_0^T (\vph_x+\psi+l\omega)^2(0)  dt+\frac{k_2|g(L)|}{2} \int\limits_0^T (u_x+v+lw)^2(L)  dt+
		\left(\frac{k_1 g_0}{2} -\hat\varepsilon\right)\int\limits_0^T \int\limits_{\hat L}^{L_0} (\vph_x+\psi+l\omega)^2dx dt
		\\+\left(\frac{k_2 g_0}{2}-\hat\varepsilon\right) \int\limits_0^T \int\limits_{L_0}^L (u_x+v+lw)^2dx dt
		+\left(\frac{\sigma g_0}{2} -\hat\varepsilon\right)\int\limits_0^T \int\limits_{\hat L}^{L_0}(\omega_x-l\vph)^2 dx dt\\+\left(\frac{\sigma g_0}{2} -\hat\varepsilon\right)\int\limits_0^T \int\limits_{L_0}^L(w_x-lu)^2 dx dt+\left(\frac{\sigma g(0)}{2}-\hat \varepsilon\right) \int\limits_0^T (\omega_x-l\vph)^2(0)  dt\\+\frac{\sigma |g(L)|}{2} \int\limits_0^T (w_x-lu)^2(L)  dt		
		+\left(\frac{\nu_1g_0}{2}-\hat\varepsilon\right)\int\limits_0^T \int\limits_{\hat L}^{L_0} \psi_x^2  dx dt\\+\left(\frac{\nu_2g_0}{2}-\hat\varepsilon\right)\int\limits_0^T \int\limits_{L_0}^L v_x^2  dx dt+\frac{\nu_1(\nu_2-\nu_1)|g(L_0)|}{2\nu_2}\int\limits_0^T  \psi_x^2(L_0) dt+\left(\frac{\nu_1g(0)}{2}-\hat\varepsilon\right)\int\limits_0^T  \psi_x^2(0) dt+
		\frac{\nu_2|g(L)|}{2}\int\limits_0^T  v_x^2(L) dt
		\\\le 
		\left(\frac{\nu_1g_1}{2}+\hat\varepsilon\right)\int\limits_0^T \int\limits_0^{\hat L} \psi_x^2  dx dt+\left(\frac{\sigma g_1}{2}+\hat\varepsilon\right) \int\limits_0^T \int\limits_0^{\hat L}(\omega_x-l\vph)^2 dx dt+\left(\frac{k_1 g_1}{2} +\hat\varepsilon\right)\int\limits_0^T \int\limits_0^{\hat L} (\vph_x+\psi+l\omega)^2dx dt\\+	\frac{\rho_1 g_1}{2} \int\limits_0^T \int\limits_0^{\hat L}\vph_t^2 dx dt+\frac{\beta_1 g_1}{2} \int\limits_0^T\int\limits_0^{\hat L}\psi_t^2 dx dt+\frac{\rho_1g_1}{2}  \int\limits_0^T\int\limits_0^{\hat L} \omega_t^2 dx dt\\+ C(E(0)+E(T))+C\int\limits_0^T \int\limits_{L_0}^L(\theta_x^2+\xi_x^2) dxdt+C(R,T)lot.
	\end{multline}
	{\it Step 2.} 
	Now we define $z$ and $p$ as solutions to the elliptic problems
	\begin{gather}
		-z_{xx}+z=\psi_t, \quad x\in (0,L_0),\label{z1}\\
		z_x(0)=0,\quad z_x(L_0)=0 \label{z2}
	\end{gather}
	and
	\begin{gather}
		-p_{xx}+p=\theta , \quad x\in (0,L_0), \label{p1}\\
		p(0)=0,\quad p(L_0)=0. \label{p2}
	\end{gather}
	Multiplying  equation \eqref{eq32} by $z_x$, integrating over $(0,L_0)$ and $(0,T)$, integrating by parts and using \eqref{z1} we obtain
	\begin{multline}
		\label{est1}
		-\delta \int\limits_0^T \int\limits_0^{L_0}\theta  z_{tx} dx dt-\lambda \int\limits_0^T \int\limits_0^{L_0}\theta_x  \psi_t   dx dt+\lambda \int\limits_0^T \int\limits_0^{L_0}\theta_x  z  dx dt-\alpha_2 \int\limits_0^T \int\limits_0^{L_0}\psi_t z dx dt\\+ \alpha_2 \int\limits_0^T \int\limits_0^{L_0}\psi_t^2 dx dt= \delta  \int\limits_0^{L_0}\theta(x,0)  z_{x}(x,0) dx-\delta  \int\limits_0^{L_0}\theta(x,T)  z_{x}(x,T) dx 
	\end{multline}
	Subsituting \eqref{p1} into the first term in \eqref{est1} and integrating by parts we obtain
	\begin{multline}
		\label{est11}
		\delta \int\limits_0^T \int\limits_0^{L_0}p_x \psi_{tt} dx dt+\lambda \int\limits_0^T \int\limits_0^{L_0}\theta_x  \psi_t   dx dt+\lambda \int\limits_0^T \int\limits_0^{L_0}\theta_x  z  dx dt-\alpha_2 \int\limits_0^T \int\limits_0^{L_0}\psi_t z dx dt\\+ \alpha_2 \int\limits_0^T \int\limits_0^{L_0}\psi_t^2 dx dt= \delta  \int\limits_0^{L_0}\theta(x,0)  z_{x}(x,0) dx-\delta  \int\limits_0^{L_0}\theta(x,T)  z_{x}(x,T) dx 
	\end{multline}	
	Next we multiply \eqref{eq2} by $-\frac{\delta}{\beta_1}p_x$ and integrate by parts over $(0,L_0)$ and $(0,T)$ to get
	\begin{multline}
		\label{est2}
		-\delta\int\limits_0^T \int\limits_0^{L_0} \psi_{tt} p_x dx dt 
		+\frac{\delta\nu_1}{\beta_1}\int\limits_0^T \int\limits_0^{L_0} \psi_{x} \theta  dx dt-\frac{\delta\nu_1}{\beta_1}\int\limits_0^T \int\limits_0^{L_0} \psi_{x} p dx dt-\frac{\delta\nu_1}{\beta_1}\int\limits_0^T  \psi_{x}(0,t) p_x(0,t) dt\\+\frac{\delta\nu_1}{\beta_1}\int\limits_0^T  \psi_{x}(L_0,t) p_x(L_0,t) dt -\frac{\delta k_1}{\beta_1}\int\limits_0^T \int\limits_0^{L_0}
		(\vph_x+\psi+l\omega) p_x dx dt-\frac{\delta\alpha_1}{\beta_1}\int\limits_0^T \int\limits_0^{L_0}  \theta_xp_x dx dt\\ +\frac{\delta\alpha_1}{2\beta_1}\int\limits_0^T \int\limits_0^{L_0}\xi(\tilde\psi+\hat\psi)p_x dx dt
		-\frac{\delta\sigma}{2\beta_1} \int\limits_0^T \int\limits_0^{L_0} (\tilde\psi+\hat\psi)(\omega_x-l\varphi+\frac12(\tilde\psi+\hat\psi)\psi) p_x dx dt\\+\frac{\delta\alpha_1}{2\beta_1} \int\limits_0^T \int\limits_0^{L_0} \psi(\tilde\xi+\hat\xi)p_x  dx dt-\frac{\delta\sigma}{2\beta_1}  \int\limits_0^T \int\limits_0^{L_0} \psi(\tilde\omega_x-l\tilde\varphi+\hat\omega_x-l\hat\varphi+\frac{\hat\psi^2}{2}+\frac{\tilde\psi^2}{2}) p_xdx dt=0
	\end{multline}
	Adding the result to \eqref{est2} we obtain that for any $\hat \varepsilon>0$ 
	\begin{multline}
		\label{ol3}
		(\alpha_2-\hat \varepsilon) \int\limits_0^T \int\limits_0^{L_0}\psi_t^2 dx dt	\le C \int\limits_0^T \int\limits_0^{L_0}(\theta_x^2 +\xi_x^2+ z_x^2+z^2+p_x^2) dx dt\\+\hat \varepsilon\int\limits_0^T\int\limits_0^{L_0} ((\omega_x-l\varphi)^2+(\vph_x+\psi+l\omega)^2)dx dt+\hat \varepsilon\int\limits_0^T  (\psi_{x}^2(0,t)+\psi_{x}^2(L_0,t)) dt\\+C\int\limits_0^T  (p_x^2(0,t)+p_x^2(L_0,t))dt+
		C(E(0)+E(T))+C(R,T)lot.
	\end{multline}
	It follows from elliptic problem \eqref{p1}, \eqref{p2}  that
	\begin{equation}\label{plot} \int\limits_0^{L_0} p_x^2 dx+\int\limits_0^{L_0} p^2 dx\le C \int\limits_0^{L_0} \theta^2 dx,
	\end{equation}
	\begin{equation}\label{plot1}\int\limits_0^{L_0} \theta xp_x dx=-\int\limits_0^{L_0} p_{xx}xp_x dx+\int\limits_0^{L_0} pxp_x dx=\frac12 \int\limits_0^{L_0} p_x^2 dx-\frac12 \int\limits_0^{L_0} p^2 dx+p_x^2(L_0,t),
	\end{equation}
	and 
	\begin{equation}\label{plot2}\int\limits_0^{L_0} \theta (x-L_0) p_x dx=-\int\limits_0^{L_0} p_{xx}(x-L_0)p_x dx+\int\limits_0^{L_0} p(x-L_0)p_x dx=\frac12 \int\limits_0^{L_0} p_x^2 dx-\frac12 \int\limits_0^{L_0} p^2 dx+p_x^2(0,t).
	\end{equation}
	Therefore, \eqref{plot}--\eqref{plot2} yield
	\begin{equation}\label{plot0}
		\int\limits_0^T (p_x^2(0,t) +	p_x^2(L_0,t))dt\le C\int\limits_0^T \int\limits_0^{L_0} \theta^2 dx dt .
	\end{equation}
	Analogously, we obtain from \eqref{z1}, \eqref{z2} that for any $0<\epsilon<1/2$
	\begin{equation}\label{zlot} 	\int\limits_0^T\int\limits_0^{L_0}( z_x^2+z^2) dx dt\le C \max_{t\in[0,T]}\|\psi_t\|_{-\epsilon}^2.
	\end{equation}
	Substituting \eqref{plot2}--\eqref{zlot} into \eqref{ol3} we arrive at
	\begin{multline}
		\label{ol4}
		(\alpha_2-\hat \varepsilon) \int\limits_0^T \int\limits_0^{L_0}\psi_t^2 dx dt	\le C \int\limits_0^T \int\limits_0^{L_0}(\theta_x^2+\xi_x^2) dx dt+\hat \varepsilon\int\limits_0^T\int\limits_0^{L_0} ((\omega_x-l\varphi)^2+(\vph_x+\psi+l\omega)^2)dx dt\\+\hat \varepsilon\int\limits_0^T ( \psi_{x}^2(0,t)+\psi_{x}^2(L_0,t)) dt+
		C(E(0)+E(T))+C(R,T)lot.
	\end{multline}
	{\it Step 3.} Similarly to Step 2, we define $y$ and $q$ as solutions to elliptic problems
	\begin{gather}
		-y_{xx}+y=\omega_t, \quad x\in (0,L_0),\label{y1}\\
		y_x(0)=0,\quad y_x(L_0)=0 ,\label{y2}
	\end{gather}
	and
	\begin{gather}
		-q_{xx}+q=\xi, \quad x\in (0,L_0), \label{q1}\\
		q(0)=0,\quad q(L_0)=0. \label{q2}
	\end{gather}
	Next we multiply equation \eqref{eq31} by $y_x$, integrating by parts  with respect to $t$ and $x$ and taking into account \eqref{y1} we arrive at
	\begin{multline*}
		-\gamma \int\limits_0^T \int\limits_0^{L_0} \xi   y_{tx} dx dt+\gamma\int\limits_0^{L_0}( \xi(x,T)  y_{x}(x,T)-\xi(x,0)   y_{x}(x,0))dx+\mu \int\limits_0^T \int\limits_0^{L_0}\xi_x  y dx dt\\-\mu \int\limits_0^T \int\limits_0^{L_0}\xi_x \omega_t dx dt-l\alpha_1
		\int\limits_0^T \int\limits_0^{L_0} \varphi_t  y_x   dx dt-
		\alpha_1
		\int\limits_0^T \int\limits_0^{L_0} \omega_t  y   dx dt\\+\alpha_1
		\int\limits_0^T \int\limits_0^{L_0} \omega_t^2  dx dt+
		\frac{\alpha_1}{2}\int\limits_0^T \int\limits_0^{L_0}(\tilde\psi_t+\hat\psi_t)\psi y_x dx dt+\frac{\alpha_1}{2}\int\limits_0^T \int\limits_0^{L_0}\psi_t(\hat\psi+\tilde\psi) y_x dx dt=0.
	\end{multline*}
	Taking into account \eqref{q1} and integrating by parts we get
	\begin{multline}
		\label{ol5}
		\gamma \int\limits_0^T \int\limits_0^{L_0} q_x   \omega_{tt} dx dt+\gamma\int\limits_0^{L_0}( \xi(x,T)  y_{x}(x,T)-\xi(x,0)   y_{x}(x,0))dx+\mu \int\limits_0^T \int\limits_0^{L_0}\xi_x  y dx dt\\-\mu \int\limits_0^T \int\limits_0^{L_0}\xi_x \omega_t dx dt-l\alpha_1
		\int\limits_0^T \int\limits_0^{L_0} \varphi_t  y_x   dx dt-
		\alpha_1
		\int\limits_0^T \int\limits_0^{L_0} \omega_t  y   dx dt\\+\alpha_1
		\int\limits_0^T \int\limits_0^{L_0} \omega_t^2  dx dt+
		\frac{\alpha_1}{2}\int\limits_0^T \int\limits_0^{L_0}(\tilde\psi_t+\hat\psi_t)\psi y_x dx dt+\frac{\alpha_1}{2}\int\limits_0^T \int\limits_0^{L_0}\psi_t(\hat\psi+\tilde\psi) y_x dx dt=0.
	\end{multline}
	After multiplication of  equation  \eqref{eq3} by $-\frac{\gamma}{\rho_1}q_x$ and integration by parts with respect to $t$ and $x$ we arrive at
	\begin{multline}
		\label{o5}
		-\gamma\int\limits_0^T \int\limits_0^{L_0}\om_{tt}q_x dx dt+ \frac{\si\gamma}{\rho_1}\int\limits_0^T \int\limits_0^{L_0}(\om_x-l\vph)\xi dx dt-\frac{\si\gamma}{\rho_1}\int\limits_0^T \int\limits_0^{L_0}(\om_x-l\vph)q dx dt-\frac{\si\gamma}{\rho_1}\int\limits_0^T (\om_x-l\vph)(0,t) q_x(0,t) dt\\+\frac{\si\gamma}{\rho_1}\int\limits_0^T (\om_x-l\vph)(L_0,t) q_x(L_0,t) dt-\frac{lk_1\gamma}{\rho_1}\int\limits_0^T \int\limits_0^{L_0}(\vph_x+\psi+l\om)q_x dx dt+\frac{\sigma\gamma}{2\rho_1}\int\limits_0^T \int\limits_0^{L_0}\psi(\tilde\psi_x+\hat\psi_x)q_x dx dt\\+\frac{\sigma\gamma}{2\rho_1}\int\limits_0^T \int\limits_0^{L_0}(\tilde\psi+\hat\psi)\psi_xq_x dx dt-\frac{\alpha_1\gamma}{\rho_1}\int\limits_0^T \int\limits_0^{L_0}\xi_x q_x dx dt=0.
	\end{multline}
	Adding \eqref{o5} to \eqref{ol5} and using the fact that
	analogously to \eqref{plot}, \eqref{plot0} and \eqref{zlot} we have
	$$ \int\limits_0^T\int\limits_0^{L_0}( q_x^2 + q^2 )dx dt+\int\limits_0^T (q_x^2(0,t) +	q_x^2(L_0,t))dt\le C \int\limits_0^T\int\limits_0^{L_0} \xi^2 dx dt$$ 
	and
	$$  \int\limits_0^T\int\limits_0^{L_0}( y_x^2+y^2) dx dt\le C \max_{t\in[0,T]}\|\omega_t\|_{-\epsilon}^2,$$
	we infer
	\begin{multline*}
		(\alpha_1-\hat \varepsilon) \int\limits_0^T \int\limits_0^{L_0}\omega_t^2 dx dt	\le C \int\limits_0^T \int\limits_0^{L_0}(\theta_x^2 +\xi_x^2) dx dt+\hat \varepsilon\int\limits_0^T \int\limits_0^{L_0}(\varphi_t^2 +\psi_x^2+(\om_x-l\vph)^2+(\vph_x+\psi+l\om)^2) dx dt\\+\hat \varepsilon\int\limits_0^T  ((\om_x-l\vph)^2(0,t) +(\om_x-l\vph)^2(L_0,t))dt+
		C(E(0)+E(T))+C(R,T)lot,
	\end{multline*}
	which together with \eqref{ol4} gives
	\begin{multline}
		\label{ol6}
		(\alpha_1-\hat \varepsilon) \int\limits_0^T \int\limits_0^{L_0}\omega_t^2 dx dt+	(\alpha_2-\hat \varepsilon) \int\limits_0^T \int\limits_0^{L_0}\psi_t^2 dx dt	\le C \int\limits_0^T \int\limits_0^{L_0}(\theta_x^2 +\xi_x^2) dx dt\\+\hat \varepsilon\int\limits_0^T ( (\om_x-l\vph)^2(0,t) +\psi_{x}^2(0,t)+(\om_x-l\vph)^2(L_0,t) +\psi_{x}^2(L_0,t) )dt\\+\hat \varepsilon\int\limits_0^T \int\limits_0^{L_0}(\varphi_t^2+ \psi_x^2+(\om_x-l\vph)^2+(\vph_x+\psi+l\om)^2)dx dt+
		C(E(0)+E(T))+C(R,T)lot.
	\end{multline}
	{\it Step 4.} Now we multiply equation \eqref{eq3} by $\omega$. After integration by parts with respect to $t$ and $x$ we get
	\begin{multline*}
		-\rho_1 \int\limits_0^T \int\limits_0^{L_0}\om_{t}^2dx dt+\rho_1  \int\limits_0^{L_0}(\om_{t}(x,T)\om(x,T)-\om_{t}(x,0)\om(x,0)) dx\\+ \si \int\limits_0^T \int\limits_0^{L_0}(\om_x-l\vph)^2dx dt-\si \int\limits_0^T (\om_x-l\vph)(L_0)\omega(L_0) dt+l\si \int\limits_0^T \int\limits_0^{L_0}(\om_x-l\vph)\varphi dx dt+lk_1 \int\limits_0^T \int\limits_0^{L_0}(\vph_x+\psi+l\om)\omega dx dt\\-\frac{\sigma}{2} \int\limits_0^T \int\limits_0^{L_0}\psi(\tilde\psi_x+\hat\psi_x)\omega dx dt-\frac{\sigma}{2} \int\limits_0^T \int\limits_0^{L_0}(\tilde\psi+\hat\psi)\psi_x\omega dx dt+\alpha_1 \int\limits_0^T \int\limits_0^{L_0}\xi_x\omega dx dt=0.
	\end{multline*}		
	Therefore, multiplying this estimate by $\frac{\alpha_1}{4\rho_1}$ we have 
	\begin{multline}
		\label{ol7}
		(\frac{\alpha_1\si}{4\rho_1}-\hat \varepsilon) \int\limits_0^T \int\limits_0^{L_0}(\om_x-l\vph)^2 dx dt \le \frac{\alpha_1}{4} \int\limits_0^T \int\limits_0^{L_0}\om_{t}^2dx dt+\hat \varepsilon  \int\limits_0^T \int\limits_0^{L_0}(\psi_x^2+(\vph_x+\psi+l\om)^2) dx dt\\+\hat \varepsilon  \int\limits_0^T (\om_x-l\vph)^2(L_0) dt+C \int\limits_0^T \int\limits_0^{L_0}(\theta_x^2 +\xi_x^2) dx dt+
		C(E(0)+E(T))+C(R,T)lot.
	\end{multline}
	Analogously, we  multiply equation \eqref{eq2} by $\psi$ and integrate by parts
	\begin{multline*}
		-\be_1\int\limits_0^T \int\limits_0^{L_0}\psi_{t} ^2 dx dt+\be_1\int\limits_0^{L_0}(\psi_{t}(x,T)\psi(x,T)-\psi_{t}(x,0)\psi(x,0) )dx  \\+\nu_1 \int\limits_0^T \int\limits_0^{L_0} \psi_{x}^2 dx dt+ \nu_1 \int\limits_0^T  \psi_{x}(L_0)\psi(L_0) dx dt +k_1\int\limits_0^T \int\limits_0^{L_0}(\vph_x+\psi+l\om)\psi dx dt \\ +\alpha_2\int\limits_0^T \int\limits_0^{L_0}\theta_x\psi dx dt -\frac{\alpha_1}{2}\int\limits_0^T \int\limits_0^{L_0}\xi(\tilde\psi+\hat\psi)\psi dx dt -\frac{\alpha_1}{2}\int\limits_0^T \int\limits_0^{L_0}\psi(\tilde\xi+\hat\xi)\psi dx dt \\
		+\frac{\sigma}{2}\int\limits_0^T \int\limits_0^{L_0}(\tilde\psi+\hat\psi)(\omega_x-l\varphi+\frac12(\tilde\psi+\hat\psi)\psi)\psi dx dt +\frac{\sigma}{2}\int\limits_0^T \int\limits_0^{L_0}\psi(\tilde\omega_x-l\tilde\varphi+\hat\omega_x-l\hat\varphi+\frac{\hat\psi^2}{2}+\frac{\tilde\psi^2}{2})\psi dx dt =0.
	\end{multline*}
	Consequently,  multiplying this estimate by $\frac{\alpha_2}{4\beta_1}$
	\begin{multline}
		\label{ol8}
		(\frac{\alpha_2\nu_1}{4\beta_1}-\hat\varepsilon) \int\limits_0^T \int\limits_0^{L_0} \psi_{x}^2 dx dt\le  \frac{\alpha_2}{4}\int\limits_0^T \int\limits_0^{L_0}\psi_{t}^2 dx dt+\hat \varepsilon \int\limits_0^T \int\limits_0^{L_0}((\om_x-l\vph)^2 +(\vph_x+\psi+l\om)^2)dx dt\\+\hat \varepsilon  \int\limits_0^T \psi_x^2(L_0) dx dt+C \int\limits_0^T \int\limits_0^{L_0}(\theta_x^2 +\xi_x^2) dx dt+
		C(E(0)+E(T))+C(R,T)lot.
	\end{multline}
	Collecting  \eqref{ol6}, \eqref{ol7} and \eqref{ol8} and choosing $\hat\varepsilon$ small enough,  we obtain that for any $\bar \varepsilon>0$
	\begin{multline}
		\label{ol9}
		\int\limits_0^T \int\limits_0^{ L_0} \psi_{x}^2 dx dt+\int\limits_0^T \int\limits_0^{L_0}(\om_x-l\vph)^2 dx dt+ \int\limits_0^T \int\limits_0^{L_0}\omega_t^2 dx dt+ \int\limits_0^T \int\limits_0^{L_0}\psi_t^2 dx dt\\\le \bar \varepsilon \int\limits_0^T \int\limits_0^{L_0}(\vph_{t}^2 +(\vph_x+\psi+l\om)^2)dx dt+\bar \varepsilon\int\limits_0^T ( (\om_x-l\vph)^2(L_0,t) +\psi_{x}^2(L_0,t))dt\\+\bar \varepsilon\int\limits_0^T ( (\om_x-l\vph)^2(0,t) +\psi_{x}^2(0,t))dt+C \int\limits_0^T \int\limits_0^{L_0}(\theta_x^2 +\xi_x^2) dx dt+
		C(E(0)+E(T))+C(R,T)lot.
	\end{multline}
	{\it Step 5.} Next we multiply equation \eqref{eq2} by $\vph_x$ and integrate by parts with respect to $t$ and $x$
	\begin{multline}
		\label{ol10}
		\be_1\int\limits_0^T \int\limits_0^{L_0}\psi_{tx}\vph_{t}   dx dt-\be_1\int\limits_0^T\psi_{t}(L_0,t)\vph_{t}(L_0,t) dt+\be_1\int\limits_0^{L_0}(\psi_{t}(x,T)\vph_{x}(x,T) -\psi_{t}(x,0)\vph_{x}(x,0))dx\\+\nu_1 \int\limits_0^T \int\limits_0^{L_0}\psi_{x} \vph_{xx} dx dt-\int\limits_0^T\psi_{x}(L_0,t)\vph_{x}(L_0,t) dt+k_1\int\limits_0^T \int\limits_0^{L_0}(\vph_x+\psi+l\om)^2dx dt\\-k_1\int\limits_0^T \int\limits_0^{L_0}(\vph_x+\psi+l\om)(\psi+l\om)dx dt +\alpha_2\int\limits_0^T \int\limits_0^{L_0}\theta_x \vph_x dx dt-\frac{\alpha_1}{2}\int\limits_0^T \int\limits_0^{L_0}\xi(\tilde\psi+\hat\psi)\vph_x dx dt\\-\frac{\alpha_1}{2}\int\limits_0^T \int\limits_0^{L_0}\psi(\tilde\xi+\hat\xi)\vph_x dx dt
		+\frac{\sigma}{2}\int\limits_0^T \int\limits_0^{L_0}(\tilde\psi+\hat\psi)(\omega_x-l\varphi+\frac12(\tilde\psi+\hat\psi)\psi)\vph_x dx dt\\+\frac{\sigma}{2}\int\limits_0^T \int\limits_0^{L_0}\psi(\tilde\omega_x-l\tilde\varphi+\hat\omega_x-l\hat\varphi+\frac{\hat\psi^2}{2}+\frac{\tilde\psi^2}{2})\vph_x dx dt=0.
	\end{multline}
	Analogously, we multiply \eqref{eq1}  by $\frac{\beta_1}{\rho_1}\psi_x$ and integrate by parts
	\begin{multline}
		\label{ol11}
		-\beta_1\int\limits_0^T \int\limits_0^{L_0}\vph_{t}\psi_{tx} dx dt-\nu_1\int\limits_0^T \int\limits_0^{L_0}(\vph_x+\psi+l\om)_x \psi_x dx dt - l\si\int\limits_0^T \int\limits_0^{L_0}(\om_x-l\vph)\psi_x dx dt\\-\frac{\sigma l}{2}\int\limits_0^T \int\limits_0^{L_0}\psi(\hat\psi+\tilde\psi)\psi_x dx dt+l\alpha_1\int\limits_0^T \int\limits_0^{L_0}\xi \psi_x dx dt=0.
	\end{multline}
	Summing up \eqref{ol10} and \eqref{ol11} and using property \eqref{c1} we get the estimate
	\begin{multline}
		\label{ol12}
		\int\limits_0^T \int\limits_0^{L_0}(\vph_x+\psi+l\om)^2dx dt\le C(E(T)+E(0))+C(R,T)lot+\int\limits_0^T \int\limits_0^{L_0}(\theta_x^2+\xi_x^2)   dx dt\\\le
		C_{\varepsilon_1}\int\limits_0^T \int\limits_0^{L_0} ((\om_x-l\vph)^2+\psi_x^2) dx dt+ 	C_{\varepsilon_1}\int\limits_0^T  (\psi_x^2(L_0)+\psi_x^2(0)+\psi_t^2(L_0)) dt\\+\varepsilon_1\int\limits_0^T  ((\vph_x+\psi+l\omega)^2(L_0)+(\vph_x+\psi+l\omega)^2(0)+\vph_t^2(L_0)) dt
	\end{multline}
	for any $\varepsilon_1>0$. \\
	In what follows we denote by $C_\varepsilon$ a large constant depending on $1/\varepsilon$ for $\varepsilon$ small.
	Next we multiply \eqref{eq2} by $x\psi_x$ and  $(x-L_0)\psi_x$ successively and integtate by parts. Summing up the results, we obtain 	
	\begin{multline}
		\label{ol121}
		\int\limits_0^T(\psi_x^2(L_0)+\psi_t^2(L_0)+\psi_x^2(0)) dt\le \varepsilon_2 \int\limits_0^T \int\limits_0^{L_0}(\vph_x+\psi+l\omega)^2 dx dt+C_{\varepsilon_2} 
		\int\limits_0^T \int\limits_0^{L_0}((\om_x-l\vph)^2+\psi_x^2+\psi_t^2) dx dt\\ + C(E(T)+E(0))+C(R,T)lot+C\int\limits_0^T \int\limits_0^{L_0}(\theta_x^2+\xi_x^2) 
	\end{multline}
	for any $\varepsilon_2>0$.
	Combining estimate \eqref{ol12} with \eqref{ol121} and choosing $2\varepsilon_2<1/C_{\varepsilon_1}$ we arrive at
	\begin{multline}
		\label{ol122}
		\int\limits_0^T \int\limits_0^{L_0}(\vph_x+\psi+l\om)^2dx dt\le C(E(T)+E(0))+C(R,T)lot+C\int\limits_0^T \int\limits_0^{L_0}(\theta_x^2+\xi_x^2)   dx dt\\\le
		C_{\varepsilon_1} 
		\int\limits_0^T \int\limits_0^{L_0}((\om_x-l\vph)^2+\psi_x^2+\psi_t^2) dx dt \\+\varepsilon_1\int\limits_0^T  ((\vph_x+\psi+l\omega)^2(L_0)+(\vph_x+\psi+l\omega)^2(0)+\vph_t^2(L_0)) dt.
	\end{multline}
	\\	{\it Step 6.} Next we multiply equation \eqref{eq1} by $-\vph$ and integrate by parts to obtain
	\begin{multline*}
		\rho_1\int\limits_0^T \int\limits_0^{L_0}\vph_{t}^2 dx dt-k_1\int\limits_0^T \int\limits_0^{L_0}(\vph_x+\psi+l\om)^2 dx dt+k_1\int\limits_0^T \int\limits_0^{L_0}(\vph_x+\psi+l\om)(\psi+l\om)dx dt\\+k_1\int\limits_0^T (\vph_x+\psi+l\om)(L_0)\vph(L_0) dt + l\si\int\limits_0^T \int\limits_0^{L_0}(\om_x-l\vph)\vph dx dt\\+\frac{\sigma l}{2}\int\limits_0^T \int\limits_0^{L_0}\psi(\hat\psi+\tilde\psi)\vph dx dt-l\alpha_1\int\limits_0^T \int\limits_0^{L_0}\xi \vph dx dt=0.
	\end{multline*}
	Therefore,
	\begin{multline}
		\label{ol14}
		\frac{\rho_1}{4k_1}\int\limits_0^T \int\limits_0^{L_0}\vph_{t}^2 dx dt\le \frac12\int\limits_0^T \int\limits_0^{L_0}(\vph_x+\psi+l\om)^2 dx dt+\varepsilon_1\int\limits_0^T (\vph_x+\psi+l\om)^2(L_0) dt\\ + \varepsilon_1\int\limits_0^T \int\limits_0^{L_0}(\om_x-l\vph)^2 dx dt+ C(E(T)+E(0))+C(R,T)lot+C\int\limits_0^T \int\limits_0^{L_0}(\theta_x^2+\xi_x^2) dx dt=0.
	\end{multline}
	Combining this estimate with \eqref{ol122} we arrive at
	\begin{multline}
		\label{ol141}
		\int\limits_0^T \int\limits_0^{L_0}\vph_{t}^2 dx dt+\int\limits_0^T \int\limits_0^{L_0}(\vph_{x}+\psi+l\omega)^2 dx dt\\\le \varepsilon_3\int\limits_0^T ((\vph_x+\psi+l\om)^2(L_0)+  (\vph_x+\psi+l\om)^2(0)+\vph_t^2(L_0)) dt +C_{\varepsilon_3}\int\limits_0^T \int\limits_0^{L_0}((\om_x-l\vph)^2+\psi_x^2+\psi_t^2) dx dt\\ + C(E(T)+E(0))+C(R,T)lot+C\int\limits_0^T \int\limits_0^{L_0}(\theta_x^2+\xi_x^2) dx dt
	\end{multline}
	for any $\varepsilon_3>0$.\\
	{\it Step 8.}	
	Now we combine estimate \eqref{ol9} multiplied by $2C_{\varepsilon_3}$ and \eqref{ol141}, and choosing $2\bar\varepsilon<\varepsilon_3/C_{\varepsilon_3}$ and $\varepsilon_3<1/2$ we come to the estimate
	\begin{multline}
		\label{ol15}
		\min\{C_{\varepsilon_3},1/2\}\left(\int\limits_0^T \int\limits_0^{L_0}\vph_{t}^2 dx dt+\int\limits_0^T \int\limits_0^{L_0}(\vph_{x}+\psi+l\omega)^2 dx dt\right.\\\left.+	\int\limits_0^T \int\limits_0^{ L_0} \psi_{x}^2 dx dt+\int\limits_0^T \int\limits_0^{L_0}(\om_x-l\vph)^2 dx dt+ \int\limits_0^T \int\limits_0^{L_0}\omega_t^2 dx dt+ \int\limits_0^T \int\limits_0^{L_0}\psi_t^2 dx dt\right)\\\le \varepsilon_3\int\limits_0^T ((\vph_x+\psi+l\om)^2(L_0,t)+  (\vph_x+\psi+l\om)^2(0,t)+\vph_t^2(L_0,t)+(\om_x-l\vph)^2(L_0,t) +\psi_{x}^2(L_0,t)\\+(\om_x-l\vph)^2(0,t) +\psi_{x}^2(0,t)) dt  + C(E(T)+E(0))+C(R,T)lot+C\int\limits_0^T \int\limits_0^{L_0}(\theta_x^2+\xi_x^2) dx dt
	\end{multline}
	for any $\varepsilon_3$ small enough.
	Taking into account estimate \eqref{l01} and choosing $\varepsilon_3$ approprietly small one can infer from \eqref{ol15} that there exists $\varepsilon_0>0$ such that for any $\varepsilon<\varepsilon_0$ 
	\begin{multline}
		\label{ol16}
		\left(\int\limits_0^T \int\limits_0^{L_0}\vph_{t}^2 dx dt+\int\limits_0^T \int\limits_0^{L_0}(\vph_{x}+\psi+l\omega)^2 dx dt\right.\\\left.+	\int\limits_0^T \int\limits_0^{ L_0} \psi_{x}^2 dx dt+\int\limits_0^T \int\limits_0^{L_0}(\om_x-l\vph)^2 dx dt+ \int\limits_0^T \int\limits_0^{L_0}\omega_t^2 dx dt+ \int\limits_0^T \int\limits_0^{L_0}\psi_t^2 dx dt\right)\\\le \varepsilon\int\limits_0^T (  (\vph_x+\psi+l\om)^2(0,t) +(\om_x-l\vph)^2(0,t) +\psi_{x}^2(0,t)) dt \\ + C(E(T)+E(0))+C(R,T)lot+C\int\limits_0^T \int\limits_0^{L_0}(\theta_x^2+\xi_x^2) dx dt.
	\end{multline}
	{\it Step 9.} Consequently, it follows from \eqref{gest3} and \eqref{ol16} that  there exist  constants $M_i$, $i=\overline{\{1,3\}}$ (depending on $l$) such that
	\begin{equation}
		\label{enr1}
		\int\limits_0^T E(t) dt dt\le M_1(R) \int\limits_0^T \int\limits_0^{L_0}(\theta_x^2+\xi_x^2)dx dt\\+M_2(R,T)lot+M_3(R)(E(T)+E(0)).
	\end{equation}	
	It follows from \eqref{En1}  that there exists $C>0$ such that
	\begin{equation}
		\label{En11}
		\int\limits_0^T \int\limits_0^{L_0}(\theta_x^2+\xi_x^2) dx dt\le C\left( E(0)+ \int\limits_0^T | H(\hat U(t),\tilde U(t))|  dt\right).
	\end{equation}
	It is easy to see from the structure of $H$  that for any $\varepsilon>0$ there exists $C(\varepsilon, R,T)>0$ such that
	\begin{equation} 	\label{HEst}
		\int\limits_0^T | H(\hat U(t), \tilde U(t))|  dt\le \varepsilon\int\limits_0^T E(t)  dx dt+ C(\varepsilon, R,T)lot.
	\end{equation}
	Combining \eqref{HEst} with \eqref{En11} we arrive at
	\begin{equation}
		\label{DampEst}
		\int\limits_0^T \int\limits_0^{L_0}(\theta_x^2+\xi_x^2) dx dt\le C E(0)+C(R,T)lot+\varepsilon\int\limits_0^T E(t)  dx dt .
	\end{equation}
	Substituting \eqref{DampEst} into \eqref{enr1} we obtain
	\begin{equation}
		\label{l4}
		\int\limits_0^T E(t) dt\le  C(R,T)lot+C(R)(E(T)+E(0))
	\end{equation}	
	for some $C(R), C(R,T)>0$.\\
	Our remaining task is to estimate the last term in \eqref{En2}. It is easy to see that
	\begin{equation}\label{En22}
		\left|	\int\limits_0^T \int\limits_t^T  H(\hat U(s),\tilde U(s)) ds dt\right|\le  \int\limits_0^T E(t) dt+T^3C(R) lot.
	\end{equation}
	Then, it follows from \eqref{En2} and \eqref{En22} that
	\begin{equation}\label{En222}
		TE(T)\le  C\int\limits_0^T E(t) dt+C(T,R) lot.
	\end{equation}
	Then the combination of \eqref{En222} with \eqref{l4} leads to the estimate
	\begin{equation*}
		\label{l}
		TE(T)\le  C(R,T)lot+C(R)(E(T)+E(0)).
	\end{equation*}
	Now we fix $\varepsilon>0$. Choosing $T$ such that $T>C(R)(1+1/\varepsilon)$ one can obtain estimate \eqref{te} which together with Theorem \ref{theoremCL} immediately leads to the asymptotic smoothness of the system.
\end{proof}
\subsection{Existence of attractors.}
The following statement collects criteria on existence and properties of attractors to gradient systems.
\begin{theorem}[\cite{CFR, CL}]
	\label{abs}
	Assume that $( H, S_t)$ is a gradient asymptotically smooth  dynamical
	system. Assume its Lyapunov function $L(y)$ is bounded from above on any
	bounded subset of $H$ and the set $\cW_R=\{y: L(y) \le R\}$ is bounded for every $R$. If the
	set $\EuScript N$ of stationary points of $(H, S_t)$ is bounded, then $(S_t, H)$ possesses a compact
	global attractor. Moreover,
	the global attractor  consists of  full trajectories
	$\gamma=\{ U(t)\, :\, t\in\R\}$ such that
	\begin{equation}\label{conv-N}
		\lim_{t\to -\infty}{\rm dist}_{H}(U(t),\EuScript N)=0 ~~
		\mbox{and} ~~ \lim_{t\to +\infty}{\rm dist}_{H}(U(t),\EuScript N)=0
	\end{equation}
	and
	\begin{equation}\label{7.4.1}
		\lim_{t\to +\infty}{\rm dist}_{H}(S_tx,\EuScript N)=0
		~~\mbox{for any $x \in H$;}
	\end{equation}
	that is, any trajectory stabilizes to the set $\EuScript N$ of  stationary points.
\end{theorem}
Now we state the result on the existence of an attractor.
\begin{theorem}
	\label{th:attr}
	Let assumptions of Theorems \ref{th:grad}, \ref{th:AsSmooth},  hold true.
	Then, the dynamical system $(H, S_t)$ generated by \eqref{AEq}-\eqref{AIC} possesses a compact global attractor $\mathfrak A$ possessing properties \eqref{conv-N}, \eqref{7.4.1}.
\end{theorem}
\begin{proof}
	In view  of Theorems \ref{th:grad}, \ref{th:AsSmooth}, \ref{abs} our remaining task is to show the boundedness of the set of stationary points and the set $W_R=\{Z:  L(Z)\le R\}$, where $L$ is given by \eqref{lap}.\par
	The second statement follows immediately from the structure of the function $L$.\par
	The first statement can be easily shown by energy-like estimates  for stationary solutions.
\end{proof}
\begin{remark}
	It would be interesting to consider a more physically relevant variant of problem \eqref{Eq1}-\eqref{IC} with different paramenters in place of $\sigma$ for different parts of the beam. However, in this case one would have  nonlinear boundary conditions
	$$\si_1(\om_x-l\vph+\psi^2/2)(L_0,t)=\si_2(w_x-lu+v^2/2)(L_0,t)$$
	instead of \eqref{TC44}. It remains an open question whether it is possible to prove the existence of strong solutions or construct other smooth approximations for weak solutions in this case to perform estimates needed for the proof of asymptotic smoothness. 
\end{remark}
\begin{remark}
	Another open question is whether a compact global attractor to the dynamical system generated by  \eqref{AEq}-\eqref{AIC} exists without the assumption of the absence of external heat sources \eqref{rightside}.
\end{remark}
\section{Singular Limits on finite time intervals}
\subsection{Singular limit $l\arr 0$.}
If we formally set $l=0$ in  \eqref{AEq}-\eqref{AIC}, we obtain the contact problem for a straight Timoshenko beam
\begin{align}
	& \rho_1\vph_{tt}-k_1(\vph_x+\psi)_x=p_1(x,t),\;x\in (0,L_0), t>0,\label{Eq1l}\\
	& \be_1\psi_{tt} -\nu_1 \psi_{xx} +k_1(\vph_x+\psi) +\alpha_2\theta_x+\sigma\psi\omega_x-\alpha_1\psi\xi+\frac{\sigma}{2}\psi^3=r_1(x,t),\label{Eq2l}\\
	& \rho_1\om_{tt}- \si\om_{xx}-\sigma\psi\psi_x+\alpha_1\xi_x=q_1(x,t), \label{Eq3l}\\
	&\gamma\xi_t-\mu\xi_{xx} +\alpha_1\omega_{xt}+\alpha_1\psi_t\psi=h_1(x,t),\label{Eq31l}\\
	&\delta\theta_t-\lambda\theta_{xx}+\alpha_2\psi_{tx}=h_2(x,t)\label{Eq32l}
\end{align}
and
\begin{align}
	& \rho_2u_{tt}-k_2(u_x+v)_x =p_2(x,t), \label{Eq4l}\\
	& \be_2v_{tt}-\nu_2 v_{xx}+k_2(u_x+v)+\sigma vw_x+\frac{\sigma}{2}v^3 =r_2(x,t),\qquad x\in (L_0,L), t>0,\label{Eq5l}\\
	& \rho_2w_{tt}- \si w_{xx}-\sigma vv_x=q_2(x,t) \label{Eq6l}
\end{align}
with boundary conditions \eqref{BC}--\eqref{BC1},
trasmission conditions at point $L_0$
\begin{align}
	& \vph(L_0,t)=u(L_0,t), \quad  \psi(L_0,t)=v(L_0,t),  \quad  \om(L_0,t)=w(L_0,t),  \\
	& k_1(\vph_x+\psi)(L_0,t)=k_2(u_x+v)(L_0,t), \\
	& \nu_1 \psi_{x}(L_0,t)= \nu_2 v_{x}(L_0,t),\\
	& \om_x(L_0,t)=w_x(L_0,t),\\
	&\xi(L_0,t)=\theta(L_0,t)=0\label{TC4l}
\end{align}
and supplemented with the initial conditions.
The following theorem gives an answer, how close are solutions to \eqref{AEq}-\eqref{AIC} when $l\arr 0$, to the solution of decoupled system \eqref{Eq1l}-\eqref{TC4l}.
\begin{theorem}
	Assume that the conditions of Theorem \ref{th:WeakWP} hold.
	Let $U^{(l)}$ be the solution to \eqref{AEq}-\eqref{AIC} with  fixed $l$ and the initial data
	\begin{gather}
		\label{inl}
		U(x,0)=U_0=(\Phi_0,\Phi_1,\Theta_0)\\
		\Phi_0=(\vph_0,\psi_0,\om_0, u_0,v_0,w_0), \quad 	\Phi_1=(\vph_1,\psi_1,\om_1, u_1,v_1,w_1), \quad\Theta_0=(\xi_0,\theta_0).
	\end{gather}
	Then for every $T>0$
	\begin{align}
		&\Phi^{(l)}  \stackrel{\ast}{\rightharpoonup} \Phi=(\vph,\psi,\om, u,v,w) \quad &\mbox{in } L^\infty(0,T;H_d) \; &\mbox{ as } l\arr 0,\label{scl1}\\
		&\Phi^{(l)}_t  \stackrel{\ast}{\rightharpoonup} \Phi_t=(\vph_t,\psi_t,\om_t, u_t,v_t,w_t) \quad &\mbox{in } L^\infty(0,T;H_v)\; &\mbox{ as } l\arr 0,\\
		&\Theta^{(l)}  \stackrel{\ast}{\rightharpoonup} \Theta=(\xi,\theta) \quad &\mbox{in } L^\infty(0,T;H_\theta) \; &\mbox{ as } l\arr 0\\
		&\Phi^{(l)}  \to \Phi=(\vph,\psi,\om, u,v,w) \quad &\mbox{in } C(0,T;H_v) \; &\mbox{ as } l\arr 0.\label{scl4}
	\end{align}
	where $U=(\Phi,\Phi_t,\Theta)$ is the solution to \eqref{Eq1l}-\eqref{TC4l} with the initial conditions \eqref{inl}.
\end{theorem}

\begin{figure}[h!]
	\centering
	\includegraphics[width=0.7\textwidth, height=0.35\textheight]{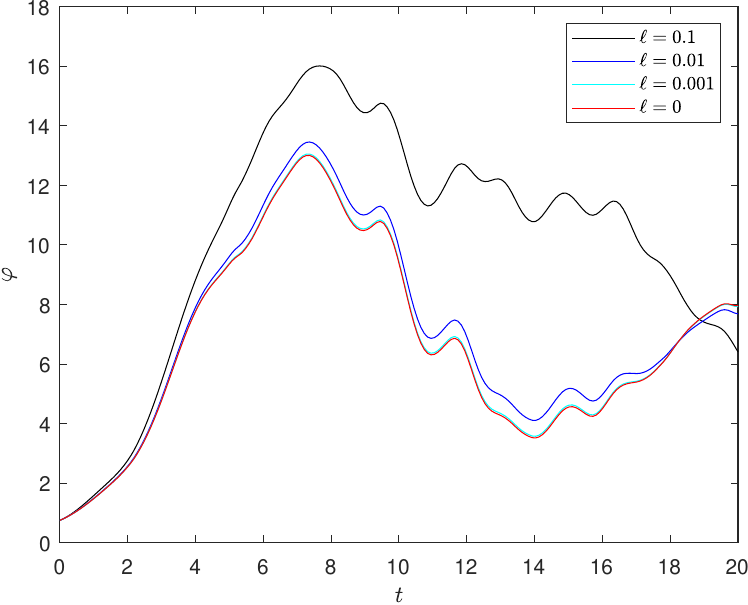}
	\caption{Transversal displacement of the beam, cross section $x=2$.}
	\label{fig:sl1_first}
\end{figure}
\begin{proof}
	The proof is similar to that of Theorem 3.1 \cite{MaMo2017} for the homogeneous Bresse beam with obvious changes, except for the limit transition in the nonlinear  term, which follows from \eqref{scl1}--\eqref{scl4}. If $l_n\to 0$ as $n\to\infty$
	\begin{multline*}
		\left|\int_{0}^T
		(\omega_{x}^{(l_n)}-l_n\varphi^{(l_n)}+\frac{1}{2}{(\psi^{(l_n)})}^2,
		\eta_{x}^1-l_n\nu^1+\psi^{(l_n)} \zeta^1)_{[0,L_0]}-(\omega_{x}+\frac12 \psi^2,\eta_{x}^1+\psi \zeta^1)_{[0,L_0]}dt
		\right|\\\le C \left|\int_{0}^T(\omega_{x}^{(l_n)}-\omega_x,\eta_x+\psi\zeta)_{[0,L_0]}dt\right|+\int_{0}^T \|\omega_{x}^{(l_n)}\|\|\psi^{(l_n)}-\psi\|\|\zeta_x\|dt\\+
		\int_{0}^T (\|\psi_x^{(l_n)}\|+\|\psi_x\|)\|\psi^{(l_n)}-\psi\|\|\eta_x+\psi\zeta\|dt\\+
		\int_{0}^T \|\psi_x^{(l_n)}\|^2\|\psi^{(l_n)}-\psi\|\|\zeta_x\|dt\\+
		|l_n|\int_{0}^T (\|\omega_{x}^{(l_n)}-l_n\varphi^{(l_n)}+\frac{1}{2}(\psi^{(l_n)})^2\|\|\nu\|+\|\varphi^{(l_n)}\|\|\eta_x\|+\|\psi^{(l_n)}\|\|\zeta_x\|)dt
		\to 0,\quad n\to \infty.
	\end{multline*} 	
\end{proof}
\begin{figure}[h!]
	\centering
	\includegraphics[width=0.7\textwidth, height=0.35\textheight]{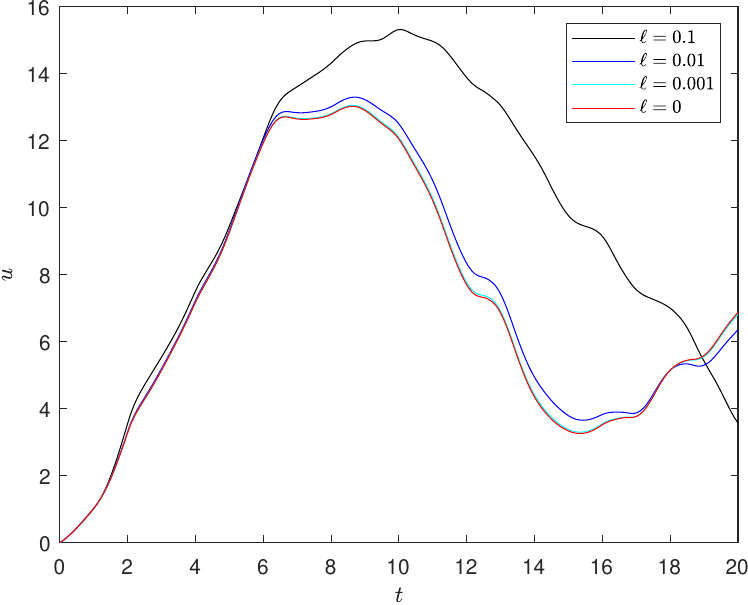}
	\caption{Transversal displacement of the beam, cross section $x=6$.}
\end{figure}
We perform numerical modelling for the original problem with $l= 1/10,  1/100, 1/1000$ and the limiting problem ($l=0$) with the following values of constants $\rho_1=\rho_2=\alpha_1=\alpha_2=1$, $\beta_1=\beta_2=2$,  $\sigma=2$, $\nu_1=4$, $\nu_2=8$, $k_1=1$, $k_2=4$, $\delta=\mu=\gamma=\lambda=1$, $L=10$, $L_0=4$ and the right-hand side
\begin{align}
	& p_1(x)=\sin x, && r_1(x)=x, && q_1(x)=\sin x, && h(x)= x, \label{rhs1}\\
	& p_2(x)=\cos x, && r_2(x)=x+1, && q_2(x)=\cos x, && g(x)=0 \label{rhs2}
\end{align}
\begin{figure}[h!]
	\centering
	\includegraphics[width=0.7\textwidth, height=0.35\textheight]{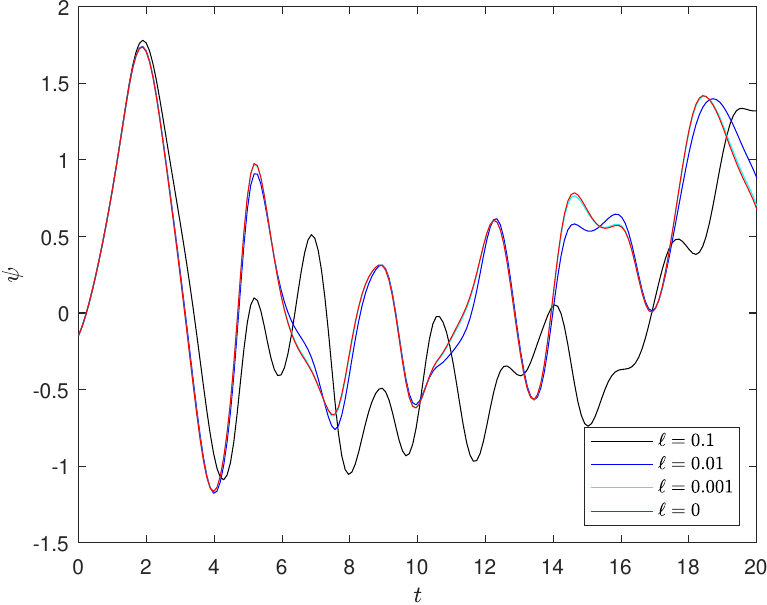}
	\caption{Shear angle variation of the beam, cross section $x=2$.}
\end{figure}
\begin{figure}[h!]
	\centering
	\includegraphics[width=0.7\textwidth, height=0.35\textheight]{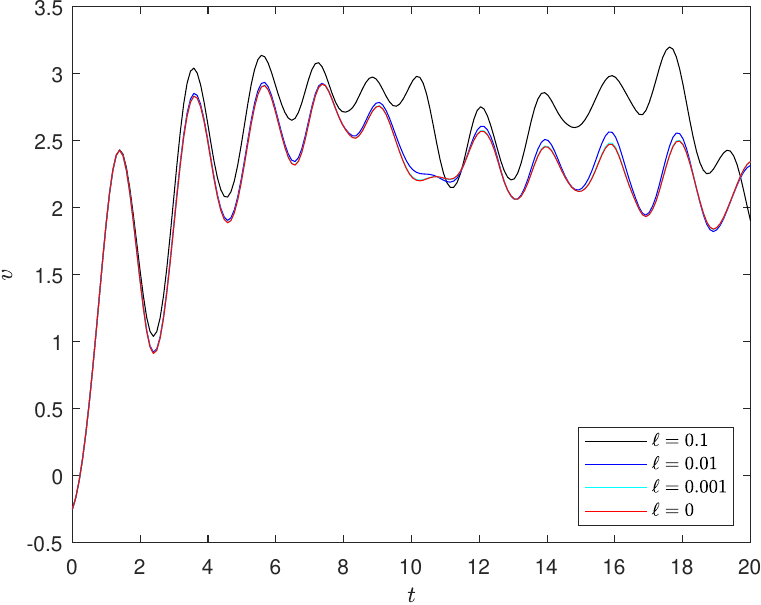}
	\caption{Shear angle variation of the beam, cross section $x=6$.}
\end{figure}
\begin{figure}[h!]
	\centering
	\includegraphics[width=0.7\textwidth, height=0.35\textheight]{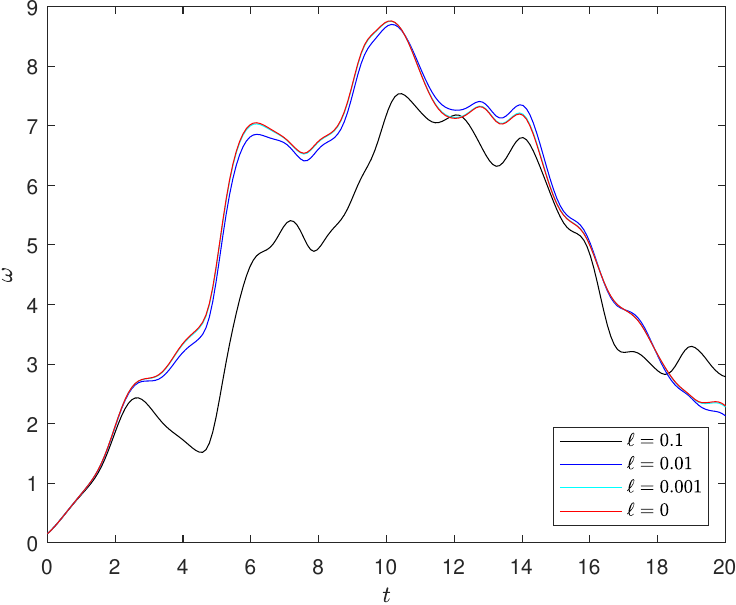}
	\caption{Longitudinal displacement of the beam, cross section $x=2$.}
\end{figure}
\begin{figure}[h!]
	\centering
	\includegraphics[width=0.7\textwidth, height=0.35\textheight]{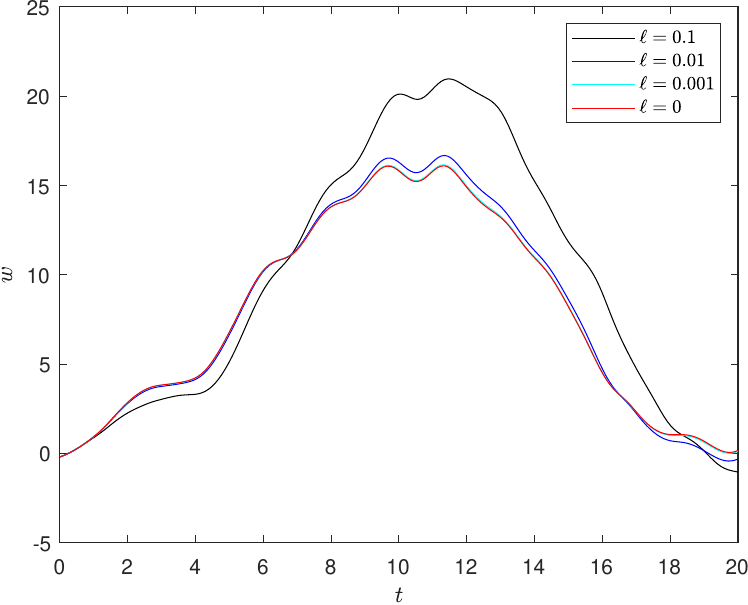}
	\caption{Longitudinal displacement of the beam, cross section $x=6$.}
	\label{fig:sl1_last}
\end{figure}
and the following initial data:
\begin{align*}
	&\vph_0(x)=-\frac{3}{16}x^2+\frac{3}{4}x,&& \vph_1(x)=\frac{x}{4}, && u_0(x)=0, && u_1(x)=-\frac{1}{6}(x-10),\\
	&\psi_0(x)=\frac{1}{192}x^2-\frac{1}{12}x,&& \psi_1(x)=\frac{x}{4}, && v_0(x)=\frac{1}{96}x^2-\frac{5}{48}x,&& v_1(x)=-\frac{1}{6}(x-10), \\
	&\om_0(x)=-\frac{3}{80}x^2+\frac{3}{20}x,&& \om_1(x)=\frac{x}{4}, && w_0(x)=\frac{1}{40}x^2-\frac {7}{20} x+1,&& w_1(x)=-\frac{1}{6}(x-10),\\
	& \xi_0(x)=x^2-4x, && \theta_0(x)=x^2-4x.
\end{align*}
Figures \ref{fig:sl1_first}-\ref{fig:sl1_last} show the behavior of solutions when $l\arr 0$ for the chosen cross sections of the beam.

\subsection{Singular limit $k_i\arr \infty, \;l\arr 0$.}
The singular limit for a straight Timoshenko beam system ($l=0$)  as $k_i\arr +\infty$ is the Euler-Bernoulli beam equation \cite[Ch. 4]{Lag1989}. We have a similar result for system \eqref{Eq1}-\eqref{IC} in case of the simultaneous limit transition $k_i\arr \infty, \;l\arr 0$. In this case the singular limit is a full von Kármán system with transmission boundary conditions (see, e.g. \cite{Fast2023}). 

\begin{theorem}
	Let  the conditions of Theorem \ref{th:WeakWP} be satisfied and, moreover,
	\begin{align}
		\begin{split}
			(\vph_0,u_0)\in\left\{\vph_0\in H^2(0,L_0), \; u_0\in H^2(L_0,L), \; \vph_0(0)=u_0(L)=0,\;\right. \\
			\left. \partial_x\phi_0(0)=\partial_x u_0(L)=0,\; \partial_x\vph_0(L_0,t)=\partial_x u_0(L_0,t) \right\};
		\end{split}\label{ICDsmooth}\\
		& \psi_0=-\pd_x\vph_{0},\;v_0=-\pd_x u_{0}; \label{ICdep} \\
		& (\vph_1,u_1) \in\{\vph_1\in H^1(0,L_0), \; u_1\in H^1(L_0,L),\; \vph_1(0)=u_1(L)=0,\; \vph_1(L_0,t)= u_1(L_0,t)\}; \label{ICVSmooth}  \\
		& \omega_0=w_0=0; \label{ICLZero}  \\
		&	r_1\in L^\infty(0,T;H^1(0,L_0)), \; r_2\in L^\infty(0,T;H^1(L_0,L)), \;
		r_1(L_0,t)=r_2(L_0,t) \; \mbox{ for a.a. } t>0.
		\label{Rsmooth}
	\end{align}
	Let  $k_j^{(n)}\arr \infty$, $l^{(n)}\arr 0$ as $n\arr \infty$, and $U^{(n)}=(\Phi^{(n)},\Phi_t^{(n)},\Theta^{(n)})$ be the solutions to \eqref{AEq}-\eqref{AIC} with fixed $k_j^{(n)}, \; l^{(n)}$ and the same initial data $U_0=(\Phi_0,\Phi_1,\Theta_0)$
	\begin{align*}
		&\Phi(x,0)=\Phi_0=(\vph_0,\psi_0,\om_0, u_0,v_0,w_0)(x), \\	&\Phi_t(x,0)=\Phi_1=(\vph_1,\psi_1,\om_1, u_1,v_1,w_1), \\ &\Theta(x,0)=\Theta_0=(\xi_0,\theta_0).
	\end{align*}
	Then for every $T>0$
	\begin{align*}
		&\Phi^{(n)}  \stackrel{\ast}{\rightharpoonup} (\vph,\psi,\om, u,v,w) \quad &\mbox{in } L^\infty(0,T;H_d) \; &\mbox{ as } n\arr \infty,\\
		&\Phi^{(n)}_t  \stackrel{\ast}{\rightharpoonup} (\vph_t,\psi_t,\om_t, u_t,v_t,w_t) \quad &\mbox{in } L^\infty(0,T;H_v)\; &\mbox{ as } n\arr \infty,\\
		&\Theta^{(n)}  \stackrel{\ast}{\rightharpoonup} (\xi,\theta) \quad &\mbox{in } L^\infty(0,T;H_\theta) \; &\mbox{ as } n\arr \infty,
	\end{align*}
	where
	\begin{itemize}
		\item $\EuScript U=(W, W_t,\Theta)$, where $W= (\vph, \omega, u, w )$ is the solution to the thermoelastic full von Karman problem
		\begin{align}
			\begin{split}
				\rho_1\vph_{tt}-\beta_1\vph_{ttxx} +\nu_1 \vph_{xxxx}+\alpha_2\theta_{xx} -\sigma\left(\varphi_x(\omega_x+\frac{1}{2}(\varphi_x)^2)\right)_x-\alpha_1\left(\varphi_x \xi\right)_x= \qquad\\
				p_1(x,{t})+\partial_x r_1(x,{t}), \quad (x,t)\in (0,L_0)\times (0,T), \label{KirchEq1}
			\end{split}\\
			& \rho_1\om_{tt}- \si\left(\om_x+\frac{1}{2}(\varphi_x)^2\right) _x+\alpha_1\xi_x=q_1(x,t), \label{Kirchq2}\\
			&\gamma\xi_t-\mu\xi_{xx} +\alpha_1(\omega_x+\frac{1}{2}(\varphi_x)^2)_t=h_1(x,t)\label{Kirchq3}\\
			&\delta\theta_t-\lambda\theta_{xx}-\alpha_2\varphi_{txx}=h_2(x,t)\label{Kirchq4}\\
			\begin{split}
				\rho_2u_{tt}  -\beta_2u_{ttxx} +\nu_2 u_{xxxx}- \sigma\left(u_x(w_x+\frac{1}{2}(u_x)^2)\right)_x= \qquad \qquad \qquad \qquad\qquad\\
				p_2(x,{t})+\partial_x r_2(x,{t}), \quad (x,t)\in (L_0,L)\times (0,T), \label{KirchEq2}
			\end{split}\\
			& \rho_1 w_{tt}- \si\left(w_x+\frac{1}{2}(u_x)^2\right)_x=q_2(x,t), \label{Kirchq5}\\			
			&\vph(0,t)=\vph_x(0,t)=0, \, u(L,t)=u_x(L,t)=0,\\ &\om(0,t)=w(L,t)=0,\,\theta(0,t)=\xi(0,t)=0,\\
			&\vph(L_0,t)=u(L_0,t),\,   \vph_x(L_0,t)=u_x(L_0,t),\, \theta(L_0,t)=\xi(L_0,t)=0,\\
			&\nu_1 \vph_{xx}(L_0,t)=\nu_2 u_{xx}(L_0,t), 
			\label{KirchTC1}\\
			&\nu_1 \vph_{xxx}(L_0,t)-\beta_1 \vph_{ttx}(L_0,t)-\alpha_2\theta_x(L_0,t) =\nu_2 u_{xxx}(L_0,t)-\beta_2 u_{ttx}(L_0,t),\label{KirchTC2},\\
			&\omega_x(L_0,t)=w_x(L_0,t).,\label{KirchTC3}
		\end{align}
		with the initial conditions $\EuScript U=(W_0, W_1, \Theta)$
		\begin{equation*}
			W_0=(\vph_0,0, u_0,0), \; W_1=(\vph_1,\omega_1, u_0, w_0),\; \Theta_0=(\xi_0,\theta_0).
		\end{equation*}
		\item $\psi=-\vph_x,  v=-u_x$;
	\end{itemize}
\end{theorem}
\begin{figure}[h!]
	\centering
	\includegraphics[width=0.7\textwidth, height=0.35\textheight]{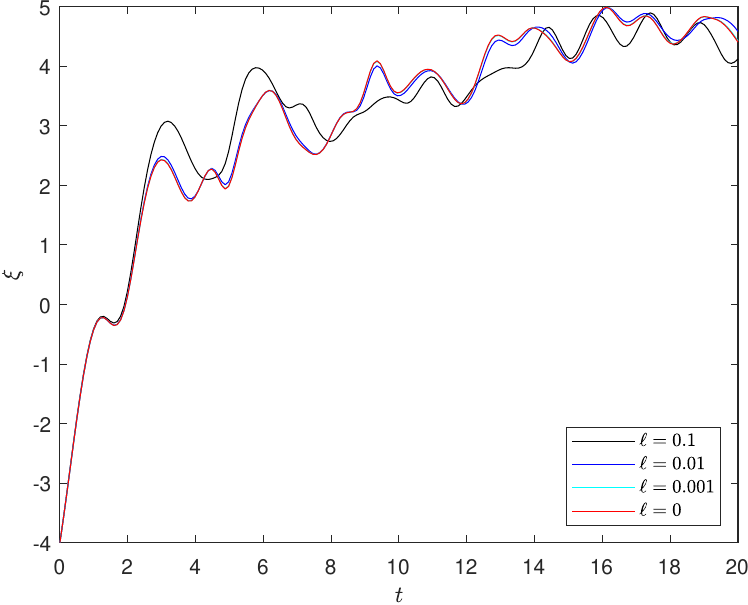}
	\caption{Temperature deviation in the longitudinal direction , cross section $x=2$.}
\end{figure}
\begin{figure}[h!]
	\centering
	\includegraphics[width=0.7\textwidth, height=0.35\textheight]{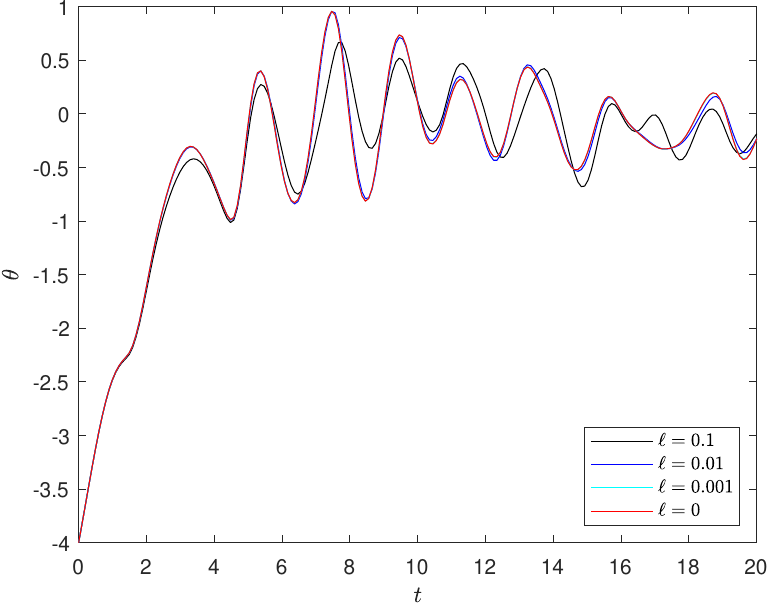}
	\caption{Temperature deviation in the vertical direction, cross section $x=2$.}
\end{figure}
\begin{proof}
	The proof uses the ideas from \cite{Fast2024, Lag1989}. We skip the details of the proof, which coincide.\\
	Energy inequality \eqref{EE} implies
	\begin{flalign}
		& \pd_t (\vph^{(n)}, \psi^{(n)}, \om^{(n)}, u^{(n)}, v^{(n)},w^{(n)})  && \mbox{ bounded in } L^\infty(0,T;H_v), \\
		& \psi^{(n)}  && \mbox{ bounded in } L^\infty(0,T;H^1(0,L_0)), \label{ConvF}\\
		& v^{(n)} && \mbox{ bounded in }  L^\infty(0,T;H^1(L_0,L)), \\
		& \om^{(n)}_x - l^{(n)}\vph^{(n)} && \mbox{ bounded in } L^\infty(0,T;L_2(0,L_0)), \\
		& w^{(n)}_x - l^{(n)} u^{(n)} && \mbox{ bounded in } L^\infty(0,T;L_2(L_0,L)), \\
		& k_1^{(n)}(\vph^{(n)}_x + \psi^{(n)} + l^{(n)} \om^{(n)}) && \mbox{ bounded in } L^\infty(0,T;L_2(0,L_0)), \\
		& k_2^{(n)}(u^{(n)}_x + v^{(n)} + l^{(n)} w^{(n)})  && \mbox{ bounded in } L^\infty(0,T;L_2(L_0,L)), \label{ConvL}\\
		& \theta^{(n)},\;\xi^{(n)}  && \mbox{ bounded in } L^\infty(0,T;L_2(0,L_0))\cap L_2(0,T;H^1(0,L_0)). \label{ConvT}
	\end{flalign}
	Thus, we can extract subsequences which converge in corresponding spaces $\ast$-weakly. Similarly to \cite{Lag1989} we have
	\begin{equation*}
		\vph^{(n)}_x + \psi^{(n)} + l^{(n)} \om^{(n)} \stackrel{\ast}{\rightharpoonup} 0 \quad \mbox{ in }  L^\infty(0,T;L_2(0,L_0)),
	\end{equation*}
	therefore
	\begin{equation*}
		\vph_x =- \psi.
	\end{equation*}
	Analogously,
	\begin{equation*}
		u_x =- v.
	\end{equation*}
	\eqref{ConvF}-\eqref{ConvL} imply
	\begin{align}
		&\om^{(n)} \stackrel{\ast}{\rightharpoonup} \om  &\mbox{ in } L^\infty(0,T;H^1(0,L_0)), &
		&w^{(n)} \stackrel{\ast}{\rightharpoonup} w  &\mbox{ in } L^\infty(0,T;H^1(L_0,L)), \label{Conv1}\\
		&\vph^{(n)} \stackrel{\ast}{\rightharpoonup} \vph  &\mbox{ in } L^\infty(0,T;H^1(0,L_0)), &
		&u^{(n)} \stackrel{\ast}{\rightharpoonup} u  &\mbox{ in } L^\infty(0,T;H^1(L_0,L)). \label{Conv2}
	\end{align}
	Thus, the Aubin's lemma gives that
	\begin{equation}\label{Conv3}
		\Phi^{(n)} \arr \Phi \mbox{ strongly in }  C(0,T; [H^{1-\ep}(0,L_0)]^3\times [H^{1-\ep}(L_0,L)]^3)
	\end{equation}
	for every $\ep>0$ and then
	\begin{equation*}
		\pd_x \vph_0 + \psi_0 + l^{(n)} \om_0 \arr 0 \quad \mbox{ strongly in }  H^{-\ep}(0,L_0).
	\end{equation*}
	This implies that
	\begin{equation*}
		\pd_x \vph_0 =- \psi_0 , \quad \om_0=0.
	\end{equation*}
	Analogously,
	\begin{equation*}
		\pd_x u_0 =- v_0 , \quad w_0=0.
	\end{equation*}
	Let us choose a test function of the form 
	$B=(\zeta_1,-\zeta_{1x},0,\zeta_2,-\zeta_{2x},0,0,0)\in F_T^0$ such that $\zeta_{1x}(L_0, t)=\zeta_{2x}(L_0, t)$ for almost all $t$. Due to \eqref{Conv1}-\eqref{Conv3}  we can pass to the limit in variational equality \eqref{VEq} as $n\arr\infty$.
	The same way as in \cite[Ch. 4.3]{Lag1989}  we obtain, that limiting functions $\vph, u$ are of higher regularity and satisfy the following variational equality
	\begin{multline}\label{LimVarEq}
		\int_0^T \int_0^{L_0}  \left(\rho_1\vph_t\zeta_{1t} - \be_1\vph_{tx}\zeta_{1tx}\right)dxdt + \int_0^T \int_{L_0}^L \left(\rho_2 u_t\zeta_{2t} - \be_1 u_{tx}\zeta_{2tx}\right)dxdt  \\
		-\int_0^{L_0}\left(\rho_1 (\vph_t \zeta_{1t})(x,0) - \be_1(\vph_{tx}\zeta_{1tx})(x,0)\right)dx + \int_{L_0}^L\left(\rho_2(u_t \zeta_{2t})(x,0) - \be_1(u_{tx}\zeta_{2tx})(x,0)\right) dx +\\
		\int_0^T \int_0^{L_0}  \nu_1\vph_{xx}\zeta_{1xx}dxdt  + \int_0^T \int_{L_0}^L\nu_2u_{xx} \zeta_{2xx}dxdt +
		\alpha_2\int_0^T\int_0^{L_0} \theta_{x}\zeta_{1x} dxdt \\
		+\sigma\int_0^T\int_0^{L_0} (\varphi_x(\omega_x+\frac12(\varphi_x)^2))\zeta_{1x}dx dt +\sigma\int_0^T\int_{L_0}^L	 (u_x(w_x+\frac12(u_x)^2))\zeta_{2x}dx dt+ \alpha_1\int_0^T\int_0^{L_0} \varphi_x\xi \zeta_{1x}dxdt  \\ =
		\int_0^T \int_0^{L_0} \left( p_1\zeta_1 - r_1 \zeta_{1x}\right)dxdt +  \int_0^T \int_{L_0}^L \left(p_2\zeta_2 - r_2 \zeta_{2x}\right)dxdt.
	\end{multline}
	Provided $\vph, u$ are smooth enough, we can integrate \eqref{LimVarEq} by parts with respect to $x,\; t$ and obtain
	\begin{multline}\label{number}
		\int_0^T \int_0^{L_0} (\rho_1\vph_{tt}-\be_1\vph_{ttxx}+\la_1\vph_{xxxx}-\alpha_2\theta_{xx}-\sigma(\varphi_x(\omega_x+\frac12(\varphi_x)^2))_x-(\varphi_x\xi)_x)\zeta_1 dxdt\\+ 
		\int_0^T \int_{L_0}^L (\rho_2 u_{tt}-\be_2\vph_{ttxx}+\la_1\vph_{xxxx}-\alpha_2\theta_{xx}-\sigma(\varphi_x(\omega_x+\frac12(\varphi_x)^2))_x-(\varphi_x\xi)_x)\zeta_2 dxdt  \\
		\int_0^T \left[\be_1\vph_{ttx}(t,L_0)-\be_2 u_{ttx}(t,L_0) -\nu_1\vph_{xxx}+\nu_2 u_{xxx}\right] (t,L_0) \be_1(t,L_0) dt  \\
		+\int_0^T  \left[\nu_1\vph_{xx}-\nu_2 u_{xx}\right] (t,L_0) \zeta_{1x}(t,L_0) dt  =
		\int_0^T \int_0^{L_0} (p_1+\pd_x r_1) \zeta_1 dxdt \\+ \int_0^T \int_{L_0}^L (p_2+\pd_x r_2) \zeta_2 dxdt + \int_0^T  \left[ r_2(t,L_0)-r_1(t,L_0)\right] \zeta_1(t,L_0)dt.
	\end{multline}
	Requiring $\zeta_1(L_0,t)$, $\zeta_{1x}(L_0,t)$ to be zero, we recover equations \eqref{KirchEq1}, \eqref{KirchEq2}  from the variational equality \eqref{number}. Then, using \eqref{Rsmooth} we get transmission conditions \eqref{KirchTC1}-\eqref{KirchTC2}. Next, we choose a test function of the form 
	$B=(0,0,0,\eta_1,0,0,\eta_2,0,0)\in F_T^0$ and again passing to the limit, and integrating by parts in \eqref{VEq}  obtain the following variational equality 
	\begin{multline}
		\int_0^T \int_0^{L_0} \left(\rho_1\omega_{tt}-\sigma (\omega_x+\frac12 (\varphi_x)^2)_x+\alpha_1\xi_x \right)\eta_{1}dx dt+\int_0^T \int_{L_0}^L \left(\rho_2w_{tt}-\sigma (w_x+\frac12 (u_x)^2)_x\right)\eta_{2}dx dt\\+
		\int_0^{T}\left[\omega_{x}(t,L_0)- w_{x}(t,L_0)\right] \eta_1(t,L_0) dt=	\int_0^T \int_0^{L_0} q_1 \eta_1 dxdt + \int_0^T \int_{L_0}^L q_2 \eta_2 dxdt.
	\end{multline}
	Therefore, choosing first $\eta_1(t,L_0)=0$ we recover equations \eqref{Kirchq2} and \eqref{Kirchq5} and subsequently transmission condition \eqref{KirchTC3}. Choosing test functions of the form 
	$B=(0,0,0,0,0,0,\tau_1,0)\in F_T^0$ and $B=(0,0,0,0,0,0,0,\tau_2)\in F_T^0$ one can also recover equations \eqref{Kirchq3} and \eqref{Kirchq4}.
\end{proof}	
We perform numerical modelling for the original problem with the initial parameters
\begin{equation*}
	l^{(1)}=0.1, \; k_1^{(1)}=0.4, \;  k_2^{(1)}=0.1;\\
\end{equation*}
We model the simultaneous convergence $l\arr 0$ and $k_1,k_2\arr \infty$ in the following way: we divide $l$ by the factor $\chi$ and multiply $k_1,k_2$ by the factor $\chi$. Calculations performed for the original problem with 
\begin{equation*}
	\chi=1,\;  \chi=10, \; \chi=100, \; \chi=1000
\end{equation*}
and the limiting problem \eqref{KirchEq1}-\eqref{KirchTC2}. 
The other constants in the original problem are the same as in the previous subsection, and we change functions in the right-hand side \eqref{rhs1}-\eqref{rhs2} as follows:
$$r_1(x)=x+4 ,  \,r_2(x)=2x.$$
We use linear dissipation $\ga(s)=s$ and we chose the following initial displacement and shear angle variation
$$
\vph_0(x)=-\frac{13}{640}x^4 +\frac{9}{40}x^3-\frac{23}{40}x^2,
$$
$$
u_0(x)=	\frac{41}{2160}x^4 -\frac{68}{135}x^3 + \frac{823}{180}x^2 -\frac{439}{27}x + \frac{520}{27}.
$$
$$
\psi_0(x)= -\left(-\frac{13}{160}x^3 +\frac{27}{40}x^2 -\frac{23}{20}x\right),
$$			 
$$
v_0(x)=-\left(\frac{41}{540}x^3 -\frac{68}{45}x^2+\frac{823}{90}x -\frac{439}{27}\right).
$$
$$\xi_0(x)= x^2-4x  ,\qquad \theta_0(x)= x^3-8x^2+16x$$
and  set 
\begin{equation*}
	\om_0(x)=w_0(x)=0.
\end{equation*}
We  choose the following initial velocities
$$
\vph_1(x)=-\frac{1}{32}x^3+\frac{3}{16}x^2, \quad u_1(x)=\frac{1}{108}x^3 -\frac{7}{36}x^2+\frac{10}{9}x-\frac{25}{27},
$$
$$\omega_1(x)=\psi_1(x)=\frac{3}{5}x,$$
$$w_1(x)=v_1(x)=-\frac{2}{5}x+4.  $$
\begin{figure}[h!]
	\centering
	\includegraphics[width=0.7\textwidth, height=0.35\textheight]{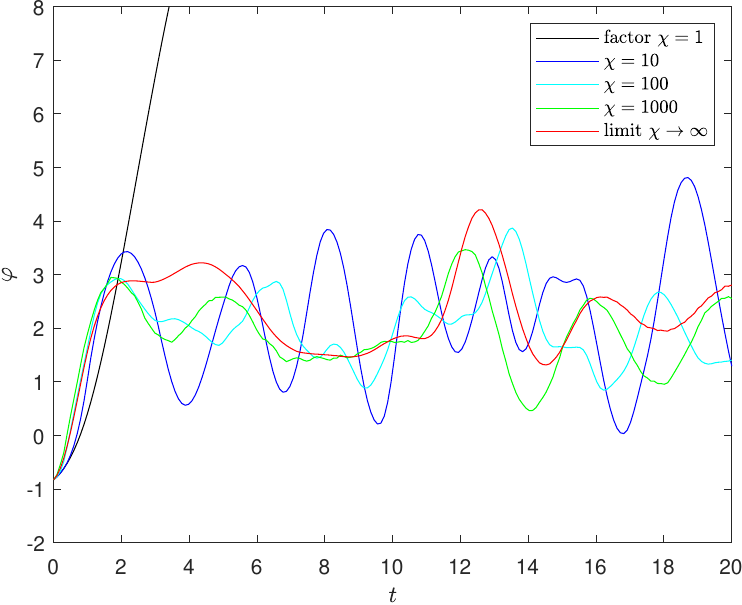}
	\caption{Transversal displacement of the beam, cross section $x=2$.}
	\label{fig:sl2_first}
\end{figure}
\begin{figure}[h!]
	\centering
	\includegraphics[width=0.7\textwidth, height=0.35\textheight]{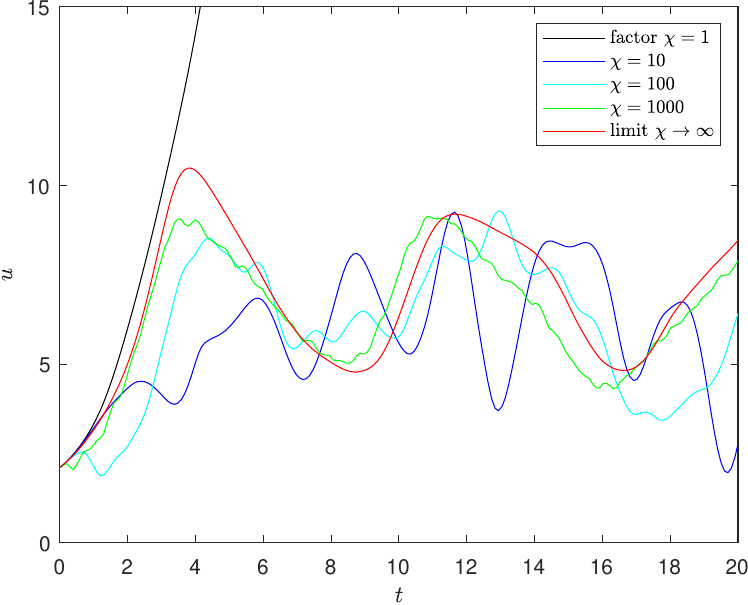}
	\caption{Transversal displacement of the beam, cross section $x=6$.}
	
\end{figure}

\begin{figure}[h!]
	\centering
	\includegraphics[width=0.7\textwidth, height=0.35\textheight]{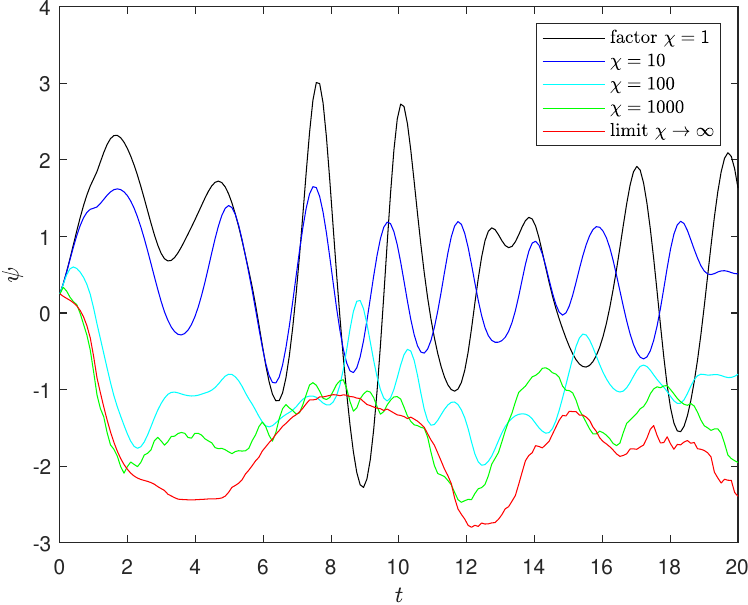}
	\caption{Shear angle variation of the beam, cross section $x=2$.}
	
\end{figure}
\begin{figure}[h!]
	\centering
	\includegraphics[width=0.7\textwidth, height=0.35\textheight]{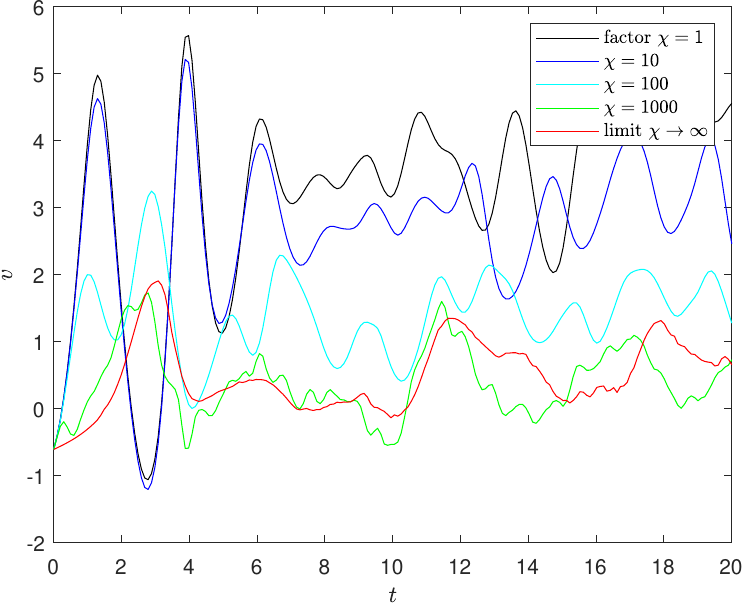}
	\caption{Shear angle variation of the beam, cross section $x=6$.}
\end{figure}

\begin{figure}[h!]
	\centering
	\includegraphics[width=0.7\textwidth, height=0.35\textheight]{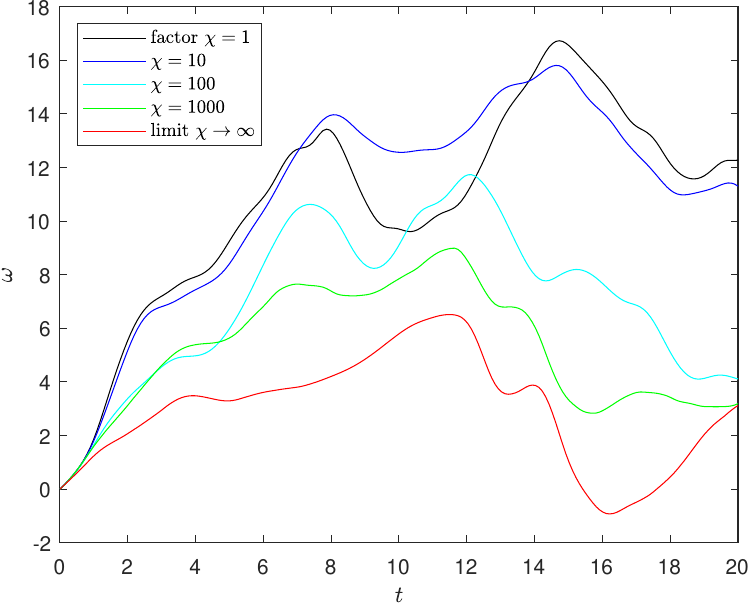}
	\caption{Longitudinal displacement of the beam, cross section $x=2$.}
\end{figure}
\begin{figure}[h!]
	\centering
	\includegraphics[width=0.7\textwidth, height=0.35\textheight]{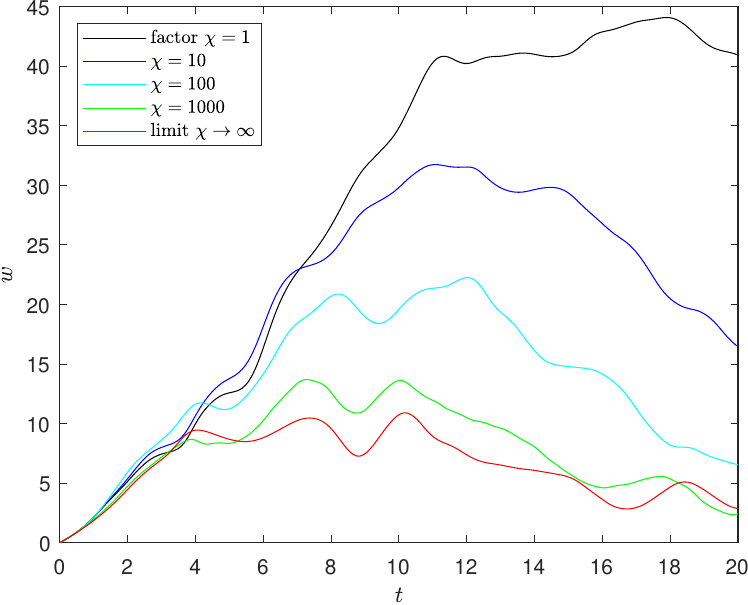}
	\caption{Longitudinal displacement of the beam, cross section $x=6$.}
	\label{fig:sl2_last}
\end{figure}
Like in \cite{Fast2024}, the double limit case appeared to be more challenging from the  point of view of numerics, than the case $l\arr 0$.
\begin{figure}[h!]
	\centering
	\includegraphics[width=0.7\textwidth, height=0.35\textheight]{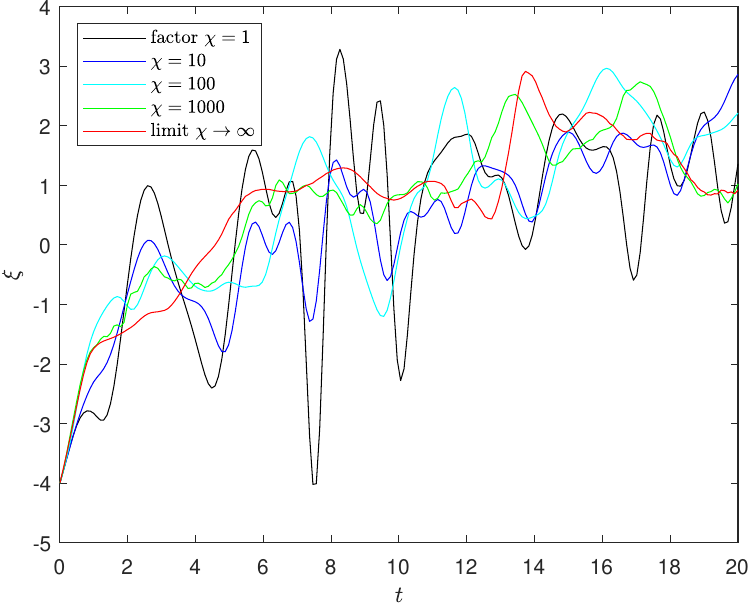}
	\caption{Temperature deviation in the longitudinal direction, cross section $x=2$.}
\end{figure}
\begin{figure}[h!]
	\centering
	\includegraphics[width=0.7\textwidth, height=0.35\textheight]{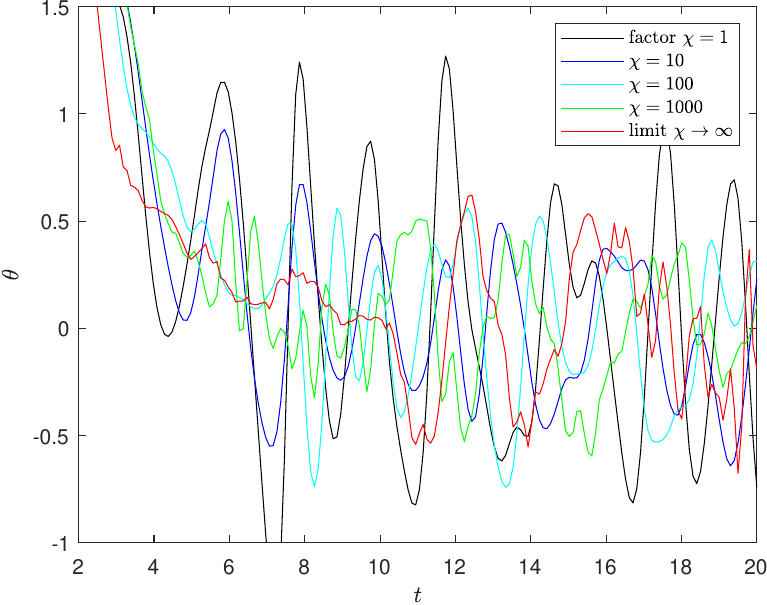}
	\caption{Temperature deviation in the vertical direction, cross section $x=2$.}
\end{figure}

Numerical simulations of the coupled system described in equations \eqref{Eq1}–\eqref{BC}, incorporating the interface conditions \eqref{TC1}–\eqref{TC4}, were carried out using a semi-discretization approach in the spatial variable $x$ for the functions  $\phi, \psi,\om, u, v,w, \xi,\theta$  and by using an
explicit scheme for the time integration. That allows to choose the discretized values at grid points near the interface in a separate step so that they obey the transmission conditions \eqref{TC1}-\eqref{TC4}. 
At each time step, only a reducible linear system involving six unknowns at three grid points (located at the interface and on either side) needed to be solved, which was computationally efficient. Any attempt to use a full implicit numerical scheme led to extremely time-expensive computations due to the large nonlinear system over all discretized values which was to solve at each time step.

Furthermore, increasing the parameters  $k_1, k_2$ raises the stiffness of the resulting system of ordinary differential equations from the semi-discretization, and the CFL-conditions requires small time steps,
otherwise numerical oscillations occur. As $k_1,\; k_2$ grow,
the beam's material becomes stiffer, and the same applies to the corresponding discretized system, although oscillations can still be observed in the resulting graphs.
The observation that the factor $\chi$ cannot be arbitrarily enlarged without challenging the numerical methods and computation times, underlines the importance
of having the limit problem for $\chi\to\infty$ in \eqref{Eq1}-\eqref{IC}.

Figures show, that the speed of convergence to the limit model in case of single limit $l\to 0$ is higher than in case of double limit  $l\to 0, k_i\to\infty$, when not only geometric configuration but also elastic properties of the beam change. 
\begin{remark}
	It is unclear whether similar singular limit results for $l\to 0$ can be obtained also for attractors of problems \eqref{Eq1}-\eqref{IC} and \eqref{Eq1l}-\eqref{inl}. 
	In the case where the beam is homogeneous and damping affects all equations in the Bresse system, this was shown in \cite{MaMo2017}. But when the damping is more localized, the current methods do not allow to prove the existence of an absorbing ball  not depending on the parameter $l$.
\end{remark}

\vskip 1000mm

{\bf Acknowlegments.}
The first author was supported by Einstein Foundation.
\vskip 1000mm

\vskip 5mm
Tamara Fastovska\\
{$^{1}$  Humboldt-Universität zu Berlin, Berlin, Germany, tamara.fastovska@hu-berlin.de \\
	$^{3}$ V.N. Karazin Kharkiv National University, Kharkiv, Ukraine}\\
\vskip 5mm
Dirk Langemann\\
$^{2}$ Technische Universität Braunschweig, Braunschweig, Germany
\end{document}